\documentclass[12pt, a4paper]{article}

\usepackage[utf8]{inputenc}
\usepackage[french, english]{babel}
\usepackage{indentfirst}
\usepackage[T1]{fontenc}
\usepackage{amsmath, amssymb, amsthm}
\usepackage{url}
\usepackage{hyperref}
\usepackage{geometry}
\usepackage{color}
\usepackage{dsfont}
\usepackage{todonotes}
\usepackage{enumitem}
\usepackage[mathscr]{euscript}
\usepackage{mathtools}
\usepackage{tikz}
\newcounter{moncompteur}

\newtheorem{theorem}[moncompteur]{Théorème}

\newtheorem{prop}[moncompteur]{Proposition}

\newtheorem{lemme}[moncompteur]{Lemme}

\renewcommand\P{\mathcal{P}}

\newcommand\R{\mathbb{R}}
\newcommand\N{\mathbb{N}}
\newcommand\C{\mathbb{C}}
\newcommand\Z{\mathbb{Z}}
\renewcommand\O[1]{O\left( #1\right)}

\renewcommand\leq{\leqslant}
\renewcommand\geq{\geqslant}
\renewcommand\theta{\vartheta}

\newcommand\re{\Re\mathfrak{e}}
\newcommand\e{\mathrm{e}}
\renewcommand\epsilon{\varepsilon}
\newcommand\chrom[1]{\MakeUppercase{\romannumeral #1}}
\renewcommand\d{\displaystyle}

\newcommand\Del{\Delta_{\bold{a,b}}}

\newcommand\m{\mathfrak{m}}

\usetikzlibrary{arrows}
\usetikzlibrary{decorations.markings}
\tikzstyle directed=[postaction={decorate,decoration={markings,
    mark=at position .035 with {\arrow{stealth}},
    mark=at position .105 with {\arrow{stealth}},
    mark=at position .25 with {\arrow{stealth}},
    mark=at position .395 with {\arrow{stealth}},
    mark=at position .537 with {\arrow{stealth}},
    mark=at position .62 with {\arrow{stealth}},
    mark=at position .765 with {\arrow{stealth}},
    mark=at position .91 with {\arrow{stealth}},
    mark=at position .98 with {\arrow{stealth}},}}]

\begin{document}
\selectlanguage{french}
\title{Lois locales de la fonction $\omega$ dans\\presque tous les petits intervalles}
\author{Élie \bsc{GOUDOUT}}
\date{}
\maketitle
\selectlanguage{english}
\vspace{-.5cm}
\begin{abstract}
For $k\geq 1$ an integer and $x\geq 1$ a real number, let $\pi_k(x)$ be the number of integers smaller than $x$ having exactly $k$ distinct prime divisors. Building on recent work of Matomäki and Radziwi\l\l, we investigate the asymptotic behavior of $\pi_k(x+h)-\pi_k(x)$ for almost all $x$, when $h$ is very small. We obtain optimal results for $k\asymp\log_2 x$ and close to optimal results for $5\leq k\leq\log_2 x$. Our method also applies to $y$-friable integers in almost all intervals $[x,x+h]$ when $\frac{\log x}{\log y}\leq (\log x)^{1/6-\varepsilon}$.\end{abstract}

\selectlanguage{french}
\section{Introduction}

\subsection{Résultats}

Soit $\mathcal{A}$ un ensemble d'entiers. On suppose que $\mathcal{A}\cap[x,2x]$ contient asymptotiquement $\delta(x)x$ entiers lorsque $x$ tend vers l'infini. Si $\mathcal{A}$ et $\delta$ ne sont pas trop erratiques et que la densité $\delta$ n'est pas trop petite, on s'attend à ce que pour tout $h\leq x$ on ait $\left\vert\mathcal{A}\cap[x,x+h]\right\vert\approx\delta(x)h$ pour presque tout $x$ dès que $\delta(x)h\rightarrow+\infty$. Récemment, Matomäki et Radziwi\l\l~\cite{matoradzi} ont développé une méthode qui permet entre autres d'obtenir un tel résultat dès que $\delta\gg 1$ et que $\mathcal{A}$ est représenté par une combinaison linéaire de fonctions multiplicatives, bornées en valeur absolue par $1$. Ils démontrent par exemple que l'ensemble des entiers ayant un nombre pair de facteurs premiers a pour densité asymptotique $\frac{1}{2}$ dans presque tous les intervalles $[x,x+h]$, lorsque $1\leq h\leq x$ tend vers l'infini. Ils montrent en fait que ce résultat est vrai dès que\footnote{On désigne par $\log_k$ la $k$-ième itérée de la fonction $\log$.} $\frac{\log_2 h}{\log h}=o\left(\delta(x)\right)$, alors que l'on espère que la condition $\delta(x)h\rightarrow+\infty$ suffise. Avec leur méthode, en utilisant l'approche de Selberg pour étudier les $\mathcal{E}_k:=\left\{n\geq 1\,:\quad\omega(n)=k\right\}$, on a ainsi montré~\cite{goudouterdoskac} que presque tout intervalle $\left[x,x+\exp\left((\log_2 x)^{1/2+\varepsilon}\right)\right]$ contenait des entiers admettant $\left\lfloor\log_2 x\right\rfloor$ facteurs premiers. Le présent article a pour but de montrer que si l'on étudie un ensemble qui vérifie certaines propriétés de crible se factorisant \og bien\fg\ --- comme par exemple les ensembles représentés par une fonction multiplicative, ou encore les $\mathcal{E}_k$ ---, alors la condition $\delta(x)h\rightarrow+\infty$, jugée optimale, suffit. L'application principale est la suivante. On note $\pi_k(x):=\left\vert\mathcal{E}_k\cap[1,x]\right\vert$ le nombre d'entiers inférieurs ou égaux à $x$ divisibles par exactement $k$ facteurs premiers distincts. Il existe une vaste littérature consacrée à l'étude des fonctions~$\pi_k$. On a notamment (\emph{cf.} Lemme~$\ref{hyp1application}$ pour un résultat plus précis), uniformément pour $x\geq 3$ et $1\leq k\ll\log_2 x$,
\[\pi_k(x)=\delta_k(x)x(1+o(1)),\]
lorsque $x$ tend vers l'infini, où l'on a posé
\begin{equation*}\label{introdefdeltak}
\delta_k(x):=\lambda\left(\kappa\right)\frac{(\log_2 x)^{k-1}}{(\log x)(k-1)!},\hspace{2cm}\kappa:=\frac{k-1}{\log_2 x},
\end{equation*}
et $\lambda(z):=\frac{1}{\Gamma(z+1)}\prod_{p\geq 2}\left(1+\frac{z}{p-1}\right)\left(1-\frac{1}{p}\right)^z$. Soient $0<r<R$ fixés. Lorsque $r<\kappa\leq R$, on pose $Q(\kappa):=\kappa\log\kappa-\kappa+1$. La formule de Stirling fournit alors
\begin{equation*}\label{eqdeltasympintro}
\delta_k(x)\asymp\frac{1}{(\log x)^{Q(\kappa)}\sqrt{\log_2 x}}.
\end{equation*}

On démontre le résultat optimal suivant.

\begin{theorem}\label{corollaire}
Soient $0<r<R$ fixés et $\psi:\R_+\rightarrow\R_+$ tendant vers l'infini en l'infini. Uniformément pour $X\geq 3$, $r\log_2 X<k\leq R\log_2 X$ et $\psi(X)\leq \delta_k(X)h\leq X$, on a
\[\pi_k(x+h)-\pi_k(x)=\delta_k(X)h\left(1+o(1)\right)\]
pour presque tout $x\sim X$ lorsque $X$ tend vers l'infini.
\end{theorem}

On note que si l'on remplace $=(1+o(1))$ par $\geq 1-\varepsilon$, on peut remplacer $\psi$ par une constante $A_{\epsilon}$ suffisamment grande. La méthode est particulièrement adaptée à l'étude des $\mathcal{E}_k$ lorsque $k\asymp\log_2 X$ puisque dans ce cas, des propriétés de crible linéaire semblables à celles de tous les entiers sont vérifiées. Lorsque $k$ est petit, ce n'est plus le cas, mais on peut tout de même partiellement adapter la méthode, et obtenir le résultat suivant. On définit, pour $k\geq 1$ et $x\geq 20$,
\begin{equation}\label{deffkdjfhkjdfghdk}
F_k(x):=\frac{(\log_2 x)^2}{k^2}\left(1-\exp\left(-\frac{k\log_3 x}{\log_2 x}\right)\right)^{-1}.
\end{equation}

\begin{theorem}\label{petitsk}
Soit $\varepsilon>0$ fixé. Uniformément pour $X\geq 20, 5\leq k\leq \log_2 X$ et $F_k(X)(\log_3 X)^{2+\varepsilon}\leq \delta_k(X)h\leq X$, on a
\[\pi_k(x+h)-\pi_k(x)\gg\frac{\delta_k(X)h}{F_k(X)}\]
pour presque tout $x\sim X$ lorsque $X$ tend vers l'infini.
\end{theorem}

Ce théorème n'est utile que pour les petites valeurs de $k$, le Théorème~$\ref{corollaire}$ fournissant un meilleur résultat lorsque $k\gg\log_2 X$. Si $k=o\left(\frac{\log_2 X}{\log_3 X}\right)$, on a $F_k(X)=\frac{(\log_2 X)^3}{k^3\log_3 X}(1+o(1))$. Quand $k$ est fixé, on perd alors essentiellement un facteur $(\log_2 X)^{3}(\log_3 X)^{1+\varepsilon}$ sur la taille minimale de $h$ et $(\log_2 X)^3(\log_3 X)^{-1}$ sur la densité espérées. Lorsque $k\geq 5$, cela précise un résultat récent de Teräväinen~\cite{teravainen}, qui a montré qu'il y avait une infinité d'entiers avec exactement $k$ facteurs premiers dans presque tous les intervalles de taille $(\log X)(\log_{k-1} X)^{C_k}$ pour $k\geq 2$ fixé et certaines constantes $C_k>0$. Pour les valeurs de $k$ supérieures à $\frac{\log_2 X}{\log_3 X}$, il est possible d'adapter la preuve du Théorème~$\ref{petitsk}$ pour améliorer légèrement le résultat. On mentionne cette possibilité à la section~\ref{kborne} sans toutefois l'approfondir dans cet article. 

Le deux théorèmes précédents sont en partie une conséquence d'un résultat d'indépendance du nombre de facteurs premiers de $a_1n+b_1$ et $a_2n+b_2$. Ce problème possède un intérêt propre. On énonce une version simplifiée du Théorème~\ref{nouvelleprop} \emph{infra}.

\begin{theorem}\label{corprinc}
Soit $R>0$ fixé. Pour tout entier $b\geq 1$ fixé, uniformément pour $x\geq 3$ et $1\leq k_1,k_2\leq R\log_2 x$, on a
\[\left\vert\left\{n\sim x\,:\quad\omega(n)=k_1,\ \omega(n+b)=k_2\right\}\right\vert\ll_b\delta_{k_1}(x)\delta_{k_2}(x)x.\]
\end{theorem}

Lorsque $k_2\asymp\log_2 x$ ou bien que $b$ est pair, ce résultat est réputé optimal à un facteur borné près. Au vu de la preuve que l'on donne, il est possible de généraliser ce résultat au cas de plusieurs translatés $n+b_1, ..., n+b_{\ell}$ pour tout $\ell\geq 3$ fixé.

Notre méthode permet de s'intéresser aux ensembles représentés par certaines fonctions multiplicatives. On traite l'exemple des entiers $x^{1/u}$-friables\footnote{On dit que l'entier $n$ est $y$-friable si son plus grand facteur premier $P^+(n)$ est inférieur ou égal à $y$.} dans de petits intervalles lorsque $u$ n'est pas trop grand. Le cas $u$ borné est directement traité par Matomäki et Radziwi\l\l\ dans~\cite{matoradzi}. On utilise les notations usuelles $\Psi(x,y):=\left\vert\left\{1\leq n\leq x\,:\quad P^+(n)\leq y\right\}\right\vert$, et $\rho$ est la fonction de Dickman, unique fonction continue sur $\R_+$ définie par $u\rho'(u)+\rho(u-1)=0$ pour $u>1$, $\rho(u)=1$ pour $0\leq u\leq 1$, et $\rho(u)=0$ pour $u<0$.

\begin{theorem}\label{thfriables}
Soient $0<\varepsilon<\frac{1}{6}$ fixé et $\psi:\R_+\rightarrow\R_+$ tendant vers l'infini en l'infini. Alors uniformément pour $X\geq 3, 1\leq u\leq (\log X)^{1/6-\varepsilon}$ et $\left(1+\rho(u)^{-1}\right)^{\psi(X)}\leq h\leq X$, on a
\[\Psi\left(x+h,X^{1/u}\right)-\Psi\left(x,X^{1/u}\right)=\rho(u)h(1+o(1))\]
pour presque tout $x\sim X$ lorsque $X$ tend vers l'infini.
\end{theorem}

De même que pour le Théorème~$\ref{corollaire}$, si on remplace $=(1+o(1))$ par $\geq 1-\varepsilon$, on peut remplacer $\psi$ par une constante $A_{\epsilon}$ suffisamment grande. Ce théorème est \emph{a priori} \og plus faible\fg\ (c'est-à-dire plus loin de ce que l'on espère) que le Théorème~$\ref{corollaire}$. Cela tient à la difficulté, à l'heure actuelle, de majorer efficacement un cardinal du type $\left\vert\left\{n\sim x\,:\quad P^+(n(n+1))\leq x^{1/u}\right\}\right\vert$.

Pour les Théorèmes~\ref{corollaire},~\ref{petitsk} et~\ref{thfriables}, il est possible d'obtenir une borne pour le cardinal de l'ensemble exceptionnel. Cette borne n'est pas très bonne en général, et on ne s'y est pas intéressé. Hildebrand et Tenenbaum~\cite[Theorem~$5.7$]{HildTenensurvolfriables} ont montré que lorsque $h\geq x^{1/u}\exp\left((\log X)^{1/6}\right)$, le Théorème~$\ref{thfriables}$ est valable avec un ensemble exceptionnel de mesure $\ll_{\varepsilon}X\exp\left((\log X)^{1/6-\varepsilon}\right)$. Ils obtiennent par ailleurs un terme d'erreur pour l'estimation asymptotique.

Enfin, on obtient le résultat suivant sur tous les intervalles, qui est une extension au cas $u$ non borné d'un théorème de Matomäki et Radziwi\l\l~\cite{matoradzi}, utilisant les mêmes outils.

\begin{theorem}\label{tsintervallesfriables}
Soit $0<\varepsilon<\frac{1}{6}$ fixé. Il existe une constante $u_0=u_0(\varepsilon)$ telle que pour $x\geq 1$, $u_0\leq u\leq (\log x)^{1/6-\varepsilon}$ et $\rho(u)^{-3-\varepsilon}\leq h\leq\sqrt{x}$ on ait
\[\Psi\left(x+h\sqrt{x},x^{1/u}\right)-\Psi\left(x,x^{1/u}\right)\geq\rho\left(u\right)^2\frac{h\sqrt{x}}{(\log x)^3}.\]
\end{theorem}

Pour $1\leq u\leq u_0$, cela implique que l'on a $\Psi\left(x+h\sqrt{x},x^{1/u}\right)-\Psi\left(x,x^{1/u}\right)\geq\rho\left(u_0\right)^2\frac{h\sqrt{x}}{(\log x)^3}$ pour $\rho(u_0)^{-3-\varepsilon}\leq h\leq\sqrt{x}$. L'exposant $3+\varepsilon$ peut être remplacé par $\frac{5c+1}{2}+\varepsilon$, si $c>0$ est une constante pour laquelle on sait montrer certaines estimations du type $$\left\vert\left\{n\sim x\,:\quad P^+(n(n+1))\leq x^{1/u}\right\}\right\vert\ll\rho(u)^{2-c}x.$$

\subsection{Description de la méthode}

On reprend le cadre d'étude introduit par Matomäki et Radziwi\l\l\ dans~\cite{matoradzi} pour l'étude des sommes courtes de fonctions multiplicatives réelles, bornées en module par $1$. On rappelle brièvement le fonctionnement global de leur preuve, en faisant apparaître les points qui diffèrent dans notre cas. Le but est de comparer une moyenne sur l'intervalle $[x,x+h]$ à une moyenne sur $[X,2X]$. Plus précisément, il s'agit de montrer que
\begin{equation}\label{qte}
\frac{1}{X}\int_X^{2X}\left\vert\frac{1}{h}\sum_{x<n\leq x+h}f(n)-\frac{1}{X}\sum_{n\sim X}f(n)\right\vert^2\mathrm{d}x=o\left(\frac{1}{X}\sum_{n\sim X}f(n)\right)
\end{equation}
lorsque $X$ tend vers l'infini. La première étape est de travailler sur un sous-ensemble $\mathcal{S}$ de $[X,2X]$, où les entiers admettent au moins un facteur premier dans certains intervalles dépendants de $X$. Cette restriction, dans leur cas très général où aucune hypothèse de crible n'est faite, se fait au prix d'un terme d'erreur $O\left(\frac{\log_2 h}{\log h}\right)$. Ce terme d'erreur est trop grand lorsque $\sum_{n\sim X}f(n)=o(X)$, par exemple lorsque $f$ est l'indicatrice d'un ensemble de densité $\delta$ tendant vers $0$. Cependant, si cet ensemble vérifie certaines propriétés de crible linéaire, on peut essentiellement se ramener à un $O\left(\delta\frac{\log_2 h}{\log h}\right)=o(\delta)$, lorsque $h$ tend vers l'infini. L'étape suivante consiste à relier la quantité~$(\ref{qte})$ à une intégrale du polynôme de Dirichlet associé à $f$ sur $\mathcal{S}\cap[X,2X]$. Cette étape ne pose aucun problème dans notre cas. La majoration de la dernière intégrale repose alors en grande partie sur la décomposition du polynôme en produit de deux autres polynômes, l'un de petite longueur, ayant un support inclus dans l'ensemble des nombres premiers, et le reste, de taille proche de celle du polynôme initial. Cette factorisation est aisée grâce à la structure de $\mathcal{S}$, et la multiplicativité de $f$. Dans le cas où $f$ est l'indicatrice d'un ensemble, elle n'est pas toujours multiplicative. On impose alors à l'ensemble de vérifier une propriété de factorisation pour procéder de la même manière. Ainsi, on peut par exemple étudier les ensembles $\mathcal{E}_k$. La dernière étape consiste essentiellement à montrer que le polynôme court prend de petites valeurs --- on peut alors le sortir de l'intégrale ---, puis à majorer l'intégrale du polynôme restant. Le polynôme court apportant déjà un facteur $o(1)$, il suffit de majorer celle-ci avec la borne \og triviale\fg\ pour ce problème, qui est fournie par le Lemme~$\ref{MVT}$. On remarque cependant que ce faisant, on perd un facteur $\delta^{-1}$, ce qui est gênant lorsque $\delta\rightarrow 0$. Pour rétablir cette perte, on emploie alors le Lemme~$\ref{IMVT}$, qui est une amélioration du Lemme~$\ref{MVT}$. On peut alors recouvrir le facteur $\delta^{-1}$, à condition que l'ensemble sur lequel on somme vérifie non seulement des conditions de crible linéaire, mais aussi en dimension $2$. 

Lorsqu'il s'agit de montrer le Théorème~$\ref{petitsk}$, on est amené à considérer un cadre légèrement différent. En effet, les entiers qui ont peu de facteurs premiers ne vérifient pas les mêmes propriétés de crible que des entiers génériques.

Dans la section~\ref{lemmesutiles}, on énonce plusieurs lemmes utiles. À la section~\ref{section3}, on définit le cadre adapté et on énonce le Théorème~$\ref{th1}$, avant de le démontrer, dans la section~\ref{sectionth1}. Les sections~\ref{applis} et~\ref{applifriall} sont consacrées respectivement aux Théorèmes~\ref{corollaire} et~\ref{corprinc}, et aux Théorèmes~\ref{thfriables} et~\ref{tsintervallesfriables}. Pour finir, on démontre le Théorème~$\ref{petitsk}$ à la section~\ref{kborne}.

\subsection{Remerciements}

Je tiens à exprimer toute ma gratitude à mon directeur de thèse, Régis de la Bretèche, pour son soutien précieux, ses nombreuses relectures et tous ses conseils, qui ont été essentiels à la rédaction de cet article. Je remercie par ailleurs Maksym Radziwi\l\l\ pour les échanges que nous avons eus, ainsi que Gérald Tenenbaum et Michel Balazard pour leurs remarques avisées.

\subsection{Notations}

Les lettres minuscules $p$ et $q$ sont réservées aux nombres premiers. Par conséquent, on écrit par exemple $\sum_p$ au lieu de $\sum_{p\text{ premier}}$. On désigne par $\omega$ la fonction additive qui compte le nombre de facteurs premiers distincts : $\omega(n)=\sum_{p\vert n}1$, et $\mu$ et $\varphi$ représentent respectivement les fonctions de Möbius et d'Euler. Étant donnés $E$ un ensemble d'entiers et $h$ une fonction additive, on note $h_E(n):=\sum_{\substack{p^{\nu}\Vert n \\ p\in E}}h(p^{\nu})$. Si $E$ est fini, on note $\left\vert E\right\vert$ son cardinal. On écrit $f\ll g$ ou $f=\O{g}$ (resp. $f\gg g$) pour dire qu'il existe une constante absolue $C>0$ telle que $\left\vert f\right\vert\leq C\left\vert g\right\vert$ (resp. $\left\vert f\right\vert\geq C\left\vert g\right\vert$). La région de validité de cette inégalité, si elle n'est pas précisée, est claire d'après le contexte. On écrit par exemple $\ll_{\varepsilon}$ ou $O_{\varepsilon}$ pour signifier que la constante implicite $C$ peut dépendre de $\varepsilon$. La relation $f\asymp g$ signifie que l'on a simultanément $f\ll g$ et $f\gg g$. On utilise la notation (asymétrique) $a\sim b$ pour dire $b<a\leq2b$. Pour tout entier $k\geq 2$, on désigne par $\log_k$ la $k$-ième itérée de la fonction $\log$. Enfin, pour $x>0$, on note $\log^+ x=\max\left(\log x, 0\right)$.

\section{Lemmes utiles}\label{lemmesutiles}

Les lemmes que l'on énonce ici se trouvent de manière quasiment identique dans les références localement citées. Dans cette section, pour $(a_n)_{n\geq 1}$ une suite de nombres complexes, $1\leq X\leq X'$ des réels et $s\in\C$, on note
\begin{align*}
A(s)&:=\sum_{n\sim X}\frac{a_n}{n^s},\\
A_{X,X'}(s)&:=\sum_{X<n\leq X'}\frac{a_n}{n^s}.
\end{align*}
On commence par énoncer le lemme permettant de relier la quantité qui nous intéresse initialement à une estimation sur le polynôme de Dirichlet associé (\emph{cf}.~\cite[lemma~14]{matoradzi}).

\begin{lemme}[Borne de Parseval]\label{Parseval}
Pour $X\geq 1$ et $T_0 \geq1$, on note $y_0:=\frac{X}{T_0^3}$. Pour $(a_n)_{n\geq1}$ une suite de complexes de module inférieur ou égal à $1$ et $1\leq h\leq y_0$, on a
\begin{multline*}
\frac{1}{X}\int_X^{2X}\left\vert\frac{1}{h}\sum_{x<n\leq x+h}a_n-\frac{1}{y_0}\sum_{x<n\leq x+y_0}a_n\right\vert^2\mathrm{d}x
\ll\frac{1}{T_0}+\int_{T_0}^{\frac{X}{h}}\left\vert A(1+it)\right\vert^2\mathrm{d}t\\
+\max_{T\geq\frac{X}{h}}\frac{X}{Th}\int_{T}^{2T}\left\vert A(1+it)\right\vert^2\mathrm{d}t.
\end{multline*}
\end{lemme}

On énonce maintenant le théorème de la valeur moyenne pour les polynômes de Dirichlet

\begin{lemme}\label{MVT}
Pour $(a_n)_{n\geq1}$ une suite de nombres complexes, $X\geq 1$ et $T>0$, on a
\[\int_{-T}^T\left\vert A(1+it)\right\vert^2\mathrm{d}t\ll\left(\frac{T}{X}+1\right)\frac{1}{X}\sum_{n\sim X}\left\vert a_n\right\vert^2.\]
\end{lemme}

Ce lemme, qui se trouve dans~\cite[chapter~9]{kowalskiiwaniec}, perd de son efficacité lorsque la suite $n\mapsto a_n$ est l'indicatrice d'un ensemble. Le lemme suivant, dû à Matomäki et Radziwi\l\l\ et qui se trouve dans~\cite[lemma~4]{teravainen}, corrige ce défaut.

\begin{lemme}\label{IMVT}
Pour $(a_n)_{n\geq1}$ une suite de nombres complexes, $1\leq X\leq X'$ et $T>0$, on a
\begin{equation}\label{eqIMVT}
\int_{-T}^T\left\vert A_{X,X'}(1+it)\right\vert^2\mathrm{d}t\ll T\left(\sum_{X< n\leq X'}\frac{\left\vert a_n\right\vert^2}{n^2}+\sum_{1\leq b\leq\frac{X'}{T}}\sum_{X<n,n+b\leq X'}\frac{\left\vert a_na_{n+b}\right\vert}{n(n+b)}\right).
\end{equation}
\end{lemme}

Dans le cas où $X'=2X$ et $(a_n)_n$ représente un ensemble de densité $\delta$, on s'attend en général à ce que ce lemme donne une majoration en $\left(\frac{T}{\delta }+1\right)\delta^2$, qui est exactement ce dont on aura besoin. On note que le lemme précédent aurait fournit une majoration en $\left(\frac{T}{X}+1\right)\delta$, moins forte lorsque $T\leq X$.

\begin{lemme}\label{halaszentiers}
Pour $(a_n)_{n\geq1}$ une suite de nombres complexes, $X\geq 1$ et $T>2$, si $\mathcal{T}\subset\left[-T,T\right]$ est un ensemble $1$-espacé, c'est à dire tel que $\vert t-t'\vert\geq 1$ pour tous $t\neq t'$ dans $\mathcal{T}$, alors
\[\sum_{t\in\mathcal{T}}\left\vert A(1+it)\right\vert^2\ll\left(\frac{\left\vert\mathcal{T}\right\vert\sqrt{T}}{X}+1\right)\left(\log T\right)\frac{1}{X}\sum_{n\sim X}\left\vert a_n\right\vert^2.\]
\end{lemme}

Ce lemme est utile lorsque $\left\vert\mathcal{T}\right\vert\leq T^{1/2-\epsilon}$ par exemple. Une telle information sera donnée par le lemme suivant (\emph{cf}.~\cite[lemma~8]{matoradzi}).

\begin{lemme}\label{grandesvaleurspolynome}
Pour $(a_p)_{p}$ une suite de complexes définie sur l'ensemble des nombres premiers telle que $\left\vert a_p\right\vert\leq 1$, et $P\geq 1$ un réel, on note
\[P(s):=\sum_{p\sim P}\frac{a_p}{p^s}.\]
Soit $\mathcal{T}\subset\left[-T,T\right]$ un ensemble $1$-espacé. Pour $V>0$, on désigne par $R=R\left(\mathcal{T},V\right)$ le nombre de $t\in\mathcal{T}$ tels que $\left\vert P(1+it)\right\vert\geq V^{-1}$. On a alors
\begin{equation*}
R\ll T^{2\frac{\log V}{\log P}}V^2\exp\left(2\frac{\log T}{\log P}\log_2 T\right).
\end{equation*}
\end{lemme}

Enfin, on utilise le lemme suivant, qui se trouve dans~\cite[lemma~2]{matoradziliouville} sous une forme légèrement différente, mais dont la preuve fournit en fait ce résultat.
\begin{lemme}\label{majlongpol}
Soient $\theta>\frac{2}{3}$ et $0<\varepsilon<\theta-\frac{2}{3}$ fixés. Lorsque les réels $P, Q, t$ et $X\geq 1$ vérifient $\exp\left(\left(\log X\right)^{\theta}\right)\leq P\leq Q\leq X$ et $\left\vert t\right\vert\leq X$, on a uniformément
\[\left\vert\sum_{P<p\leq Q}\frac{1}{p^{1+it}}\right\vert\ll\frac{\log X}{1+\left\vert t\right\vert}+P^{-(\log X)^{-2/3-\varepsilon}}.\]
\end{lemme}

\section{Résultat général}\label{section3}

On introduit des notations afin de définir l'ensemble $\mathcal{S}$ dont nous avons parlé dans l'introduction. Il diffère par plusieurs aspects de celui défini dans~\cite{matoradzi}. On a apporté ces modifications afin d'alléger les hypothèses du Théorème~$\ref{th1}$, en assurant une meilleure factorisation du polynôme de Dirichlet.
Pour mieux appréhender les notations à venir et le théorème suivant, il est bon d'avoir en tête les faits informels suivants. La fonction $\delta$ représente la densité de l'ensemble $\mathcal{A}$ qu'on étudie. On s'attend, de manière générale, à avoir un nombre d'entiers $n\leq X$, tels que $n\in\mathcal{A}$ et $n+1\in\mathcal{A}$, qui soit $\ll\delta^2X$. La fonction $\theta$ quantifie alors la perte que l'on a, s'il en est, sur cette dernière estimation. Pour une première approche, on peut considérer $\theta(X)=1$, ce qui revient à ne supposer aucune perte. La variable $h$ représente tout simplement la taille des intervalles dans lesquels on veut avoir des estimations asymptotiques. Pour $1\leq P\leq Q\leq X$, la proportion des entiers inférieurs à $X$ n'admettant aucun facteur premier dans l'intervalle $]P,Q]$ est $\ll\frac{\log P}{\log Q}$. On espère avoir le résultat analogue pour un ensemble d'entiers pas trop erratique. Cependant, pour les entiers inférieurs à $X$ ayant environ $\kappa\log_2 X$ facteurs premiers, la proportion est plutôt $\ll\left(\frac{\log P}{\log Q}\right)^{\kappa}$. Si on étudie les entier ayant $k\asymp\log_2 X$ facteurs premiers, il existe ainsi un $r>0$ fixé tel que l'inégalité précédente soit uniformément vérifiée en remplaçant $\kappa$ par $r$, ce qui garantit que l'ensemble $\mathcal{S}$ construit ci-dessous soit dense dans $\mathcal{A}$.
\\
\\
\indent Soient $0<r\leq 1$ et $0<\varepsilon<\frac{1}{100}$ fixés. Soient $X\geq 3$ un réel et $\delta$ et $\theta$ des fonctions vérifiant
\begin{equation}\label{defS1}
0<\delta(X)\leq1\leq\theta(X)\leq\delta(X)^{-1}.
\end{equation}
Pour $3\leq h\leq X$, on considère $\left(P_j,Q_j,H_j\right)_{j\in\N\cup\left\{\infty\right\}}$ et $K$ des réels supérieurs à $2$ vérifiant
\begin{equation}\label{defS2}
\left\{\begin{tabular}{llr}
$(\log_2 X)^2\leq K\leq\sqrt{\log X}$,\\
$\d 20\leq Q_1^{1/\log_3 Q_1}\leq P_1\leq Q_1,$&$\d Q_1=\min\left(\delta(X)\theta(X)h,\exp\left(K\right)\right)$,\\
$\d \log P_j=j^{4j/r}(\log_3 Q_1)^{2(j-1)}(\log P_1),$&$\d \log Q_j=j^{(4j+2)/r}(\log_3 Q_1)^{2(j-1)}(\log Q_1)$, &\hspace{.5cm}$\d (j\in\N)$\\
\multicolumn{2}{l}{$\d H_j=j^2\min\left(P_1^{1/6-\varepsilon}(\log Q_1)^{-1/3},\theta(X)(\log_2 Q_1)^2\right)$,}&\hspace{.5cm}$\d (j\in\N)$\\
\multicolumn{3}{l}{$\d \exp\left((\log X)^{2/3+\varepsilon}\right)\leq P_{\infty}\leq Q_{\infty}\leq\exp\left((\log X)^{1-\varepsilon}\right),$}\\
$\log P_{\infty}=o\left(\log Q_{\infty}\right)$,\\
\multicolumn{3}{l}{$\d H_{\infty}=(\log_2 X)^2$.}
\end{tabular}\right.
\end{equation}
On définit $J$ comme étant le plus grand entier $j\geq 1$ tel que $Q_j\leq\exp\left(K\right)$. On remarque que l'on a $P_1<Q_1<P_2<...<Q_J<P_{\infty}<Q_{\infty}$. On note alors
\begin{equation}\label{defS3}
\left\{\begin{tabular}{lr}
$\mathcal{I}_j:=\left[ \left\lfloor H_j\log P_j\right\rfloor,H_j\log Q_j\right],$&$\left(j\in\left[1,J\right]\cup\left\{\infty\right\}\right)$\\
$\d \mathcal{S}_j:=\bigcap_{v\in\mathcal{I}_j}\bigcap_{\substack{\e^{v/H_j}< p,q\leq\e^{(v+1)/H_j}\\P_j< p,q\leq Q_j}}\left\{n\geq 1\,:\quad \mu_{]P_j,Q_j]}^2\omega_{]P_j,Q_j]}(n)\geq 1,\ pq\nmid n\right\}.$&$\left(j\in\left[1,J\right]\cup\left\{\infty\right\}\right)$
\end{tabular}\right.
\end{equation}
Finalement, on note
\begin{equation}\label{defS4}\begin{tabular}{r}
$\d \mathcal{S}:=\bigcap_{j=1}^J\mathcal{S}_j\cap\mathcal{S}_{\infty}$.
\end{tabular}
\end{equation}

On note que l'on a $J\ll\frac{\log K}{\log_2 K}$, et donc $H_J\ll\theta(X)(\log K)^4$. Pour résumer, les entiers dans $\mathcal{S}$ ont au moins un facteur premier dans chacun des $]P_j,Q_j]$, mais pas de facteur carré dans ceux-ci. De plus, pour $j$ fixé, ils n'ont pas deux facteurs premiers trop proches dans $]P_j,Q_j]$. De ces trois propriétés, la première est la plus importante, et apparaît déjà dans~\cite{matoradzi}. Les deux autres assurent une factorisation simplifiée du polynôme de Dirichlet (\emph{cf.} Lemme~\ref{decompo}), ce qui permet d'alléger les hypothèses du théorème suivant. Enfin, on note que l'on a facilement
\[\left\vert[X,2X]\smallsetminus\mathcal{S}\right\vert\ll X\sum_{j\in\left[1,J\right]\cup\left\{\infty\right\}}\left(\frac{\log P_j}{\log Q_j}+\frac{1}{P_j}+\frac{\log\left(1+\frac{\log Q_j}{\log P_j}\right)}{H_j}\right).\]

\begin{theorem}\label{th1}
Soient $X\geq 3$ et $\mathcal{A}$ un ensemble pouvant dépendre de $X$. Soient $0<r\leq 1$ et $0<\varepsilon\leq\frac{1}{100}$ fixés. Soient $\delta$ et $\theta$ des fonctions vérifiant $(\ref{defS1})$. Pour $3\leq h\leq X$ on utilise les notations $(\ref{defS2}), (\ref{defS3})$ et la définition $(\ref{defS4})$ de $\mathcal{S}$. Soit $f:\N\rightarrow\R$ une fonction positive telle que pour $N\geq 1$ on ait $\sum_{1\leq b\leq N}f(b)\ll N$. On suppose les conditions suivantes vérifiées :
\begin{enumerate}[leftmargin=*]
\item\label{hyp1} On a $\delta(X)\gg P_{\infty}^{-(\log X)^{-2/3-\varepsilon/2}}$. Soit $T_0:=T_0(X)>1$ tel que $\frac{1}{T_0}+\frac{\theta(X)^2(\log X)^5}{T_0^2}=o\left(\delta(X)^2\right)$. On pose $y_0:=\frac{X}{T_0^3}$. Alors pour tout $\frac{X}{2}\leq x\leq 2X$ et tout $y_0\leq y\leq 2X$, on a
\begin{equation*}
\left\vert\mathcal{A}\cap\left[x,x+y\right]\right\vert=\delta(X)y\left(1+o(1)\right).
\end{equation*}
\item\label{hyp2} \emph{(crible linéaire pour $\mathcal{A}$) }
\begin{equation*}
\left\vert\left(\mathcal{A}\smallsetminus\mathcal{S}\right)\cap\left[X,2X\right]\right\vert\ll\delta(X)X\sum_{j=1}^J\left(\left(\frac{\log P_j}{\log Q_j}\right)^r+\frac{1}{P_j}+\frac{\log\left(1+\frac{\log Q_j}{\log P_j}\right)}{H_j}\right)+o\left(\delta(X)X\right).
\end{equation*}
\item\label{hyp3} \emph{(crible binaire pour $\mathcal{A}$) }Pour tout $\frac{X}{2}\leq x\leq 2X$, tout $1\leq z\leq(\log_2 X)^{2}(\log K)^{4}\theta(X)$ et tout $1\leq b\leq X^2$, on a
\begin{equation}\label{bidim1}
\left\vert\left\{n\in\mathcal{A}\,:\quad n+b\in\mathcal{A}\right\}\cap\left[x,x+\frac{x}{z}\right]\right\vert\ll f(b)\delta(X)^2\frac{x}{z}\theta(X).
\end{equation}
\item\label{hyp4} \emph{(bonne factorisation) }Pour tout $j\in\left[1,J\right]\cup\left\{\infty\right\}$ et tous $p,q\in\left]P_j,Q_j\right]$, pour tout $m\geq 1$, on a
\begin{equation*}
mq\in\mathcal{A}\cap\mathcal{S}\text{ et }p\nmid m\Longrightarrow mp\in\mathcal{A}.
\end{equation*}
\item\label{hyp5} \emph{(crible pour $\mathcal{A}'$) }On définit $\mathcal{A'}$ comme l'ensemble des quotients $m=\frac{n}{p}$ pour $n\in\mathcal{A}\,\cap\,\mathcal{S}\,\cap\,\left[X,2X\right]$ et $p\in\bigcup_{1\leq j\leq J}\left]P_j,Q_j\right]$ tels que $p\vert n$. Alors pour tout $1\leq z\leq\exp\left(K\right)$ et tout $1\leq b\leq X^2$, on a
\begin{align}
\left\vert\mathcal{A}'\cap\left[\frac{X}{z},2\frac{X}{z}\right]\right\vert&\ll\delta(X)\frac{X}{z},\nonumber\\
\left\vert\left\{n\in\mathcal{A}'\,:\quad n+b\in\mathcal{A}'\right\}\cap\left[\frac{X}{z},2\frac{X}{z}\right]\right\vert&\ll f(b)\delta(X)^2\frac{X}{z}\theta(X).\label{bidim2}
\end{align}
\item\label{hyp6} \emph{(crible pour $\mathcal{A}'$, suite) }Soit $1\leq j\leq J-1$. Soient $Y_1\in\left[P_j,Q_j\right]$ et $Y_2\in\left[P_{j+1},Q_{j+1}\right]$. On pose $\ell :=\left\lceil\frac{\log Y_2}{\log Y_1}\right\rceil$. Alors pour tous $p_1, ..., p_{\ell}\sim Y_1$, et tout $1\leq b\leq X^2$, on a
\begin{multline}\label{bidim3}
\left\vert\left\{m\sim\frac{X}{Y_2}\,:\quad m\in\mathcal{A}',\ mp_1\cdots p_{\ell}+b\in\bigcup_{q_1, ..., q_{\ell}\sim Y_1}q_1 \cdots q_{\ell}\mathcal{A}'\right\}\right\vert\\
\ll f(b)\delta(X)^2\frac{X}{Y_2}\left(\ell !\right)^{O(1)}\theta(X).
\end{multline}
\item\label{hyp7} $\theta(X)=o\left(P_1^{1/6-2\varepsilon}\right)$,
\item\label{hyp8} $\log P_1=o(\log Q_1)$.
\end{enumerate}
Alors on a
\begin{equation*}
\left\vert\mathcal{A}\cap\left[x,x+h\right]\right\vert=\delta(X)h\left(1+o(1)\right),
\end{equation*}
pour presque tout $x\sim X$, dès que $\delta(X)h$ tend vers l'infini.
\end{theorem}

\noindent\textbf{Remarques : }

\begin{enumerate}[parsep=0cm, topsep=0.2cm, label=(\roman*)]
\item Lorsque $\theta(X)=1$, la condition~\ref{hyp7} est vérifiée automatiquement. C'est ce que l'on espère en général.
\item Si $\delta(X)\gg 1$, les points \ref{hyp2}, \ref{hyp3}, \ref{hyp5} et \ref{hyp6} sont automatiquement vérifiés avec $\theta(X)=1$, et donc aussi le point~\ref{hyp7}. Il suffit donc essentiellement que $\mathcal{A}$ admette une bonne factorisation.
\item En prenant $\theta(X)=\delta(X)^{-1}$, les inégalités~$(\ref{bidim1})$,~$(\ref{bidim2})$~et~$(\ref{bidim3})$ sont automatiquement vérifiées avec $f=1$. Ce cas revient à ignorer toute la partie de crible binaire, et donc finalement à utiliser le Lemme~$\ref{MVT}$ au lieu du Lemme~$\ref{IMVT}$. On retrouve un résultat qui est à la portée de la méthode initiale de Matomäki et Radziwi\l\l.
\item Dans l'hypothèse~\ref{hyp6}, on a $\ell\ll J^{6/r}(\log_3 Q_1)^{3}$, qui est $\ll\left(\log K\right)^{O(1)}$.
\item Le point~\ref{hyp4} revient à dire que dans un intervalle $[P_j,Q_j]$, tous les facteurs premiers \og jouent le même rôle\fg. On peut alors alléger cette hypothèse en triant, dans chaque $[P_j,Q_j]$, les $p$ en selon leur rôle. Par exemple, au lieu de s'intéresser aux entiers $n$ tels que $\omega(n)=k$, on peut étudier ceux vérifiant $\tilde{\omega}(n)=\frac{3}{2}k$, où $\tilde{\omega}$ est la fonction additive valant $1$ sur les puissance des nombres premiers congrus à $1$ modulo $3$, et $2$ sur les puissances des autres nombres premiers. Notre méthode permettrait de traiter le cas des entiers représentables en somme de $2$ carrés de cette manière.
\end{enumerate}

\section{Démonstration du Théorème $\ref{th1}$}\label{sectionth1}

\subsection{Première simplification}

\begin{proof}
Soient $3\leq h\leq X$ deux réels et $\mathcal{A}$ un ensemble d'entiers. On suppose qu'il existe des paramètres $r,\varepsilon, \delta, \theta, K, P_1, P_{\infty}$ et $Q_{\infty}$ vérifiant les hypothèses du Théorème~$\ref{th1}$. Grâce à l'hypothèse \ref{hyp1}, il nous suffit de démontrer que pour presque tout $x\sim X$, on a
\begin{equation*}
\left\vert\frac{1}{h}\sum_{\substack{x<n\leq x+h\\ n\in\mathcal{A}}}1-\frac{1}{y_0}\sum_{\substack{x<n\leq x+y_0\\ n\in\mathcal{A}}}1\right\vert=o\left(\delta(X)\right)
\end{equation*}
On remarque que cela est vrai directement lorsque $h\geq y_0$ d'après l'hypothèse \ref{hyp1}. Il suffit donc de traiter le cas $h\leq y_0$. Pour cela, on travaille sur un sous-ensemble dense de $\mathcal{A}$, d'entiers qui admettent une certaine factorisation, ce qui nous permet d'utiliser la méthode de Matomäki et Radziwi\l\l. On parvient alors à démontrer le résultat suivant.

\begin{prop}\label{estimationintegrale}
Avec les notations ci-dessus, et sous les hypothèses~\ref{hyp1} à~\ref{hyp6} du Théorème~$\ref{th1}$, lorsque $h\leq y_0$ et $Q_1$ est suffisamment grand, on a
\begin{equation*}
\frac{1}{X}\int_X^{2X}\left\vert\frac{1}{h}\sum_{\substack{x<n\leq x+h\\ n\in\mathcal{A}\cap{\mathcal{S}}}}1-\frac{1}{y_0}\sum_{\substack{x<n\leq x+y_0\\ n\in\mathcal{A}\cap{\mathcal{S}}}}1\right\vert^2\mathrm{d}x\ll\delta(X)^2\left(\frac{\theta(X)}{H_1}+o(1)\right).
\end{equation*}
quand $X$ tend vers l'infini.
\end{prop}

On démontre cette proposition à la sous-section~$\ref{sssectionparseval}$. Vérifions dans un premier temps qu'elle implique bien le Théorème $\ref{th1}$. En se rappelant que $h\leq y_0$, on obtient
\begin{multline*}
\frac{1}{X}\int_X^{2X}\left\vert\frac{1}{h}\sum_{\substack{x<n\leq x+h\\ n\in\mathcal{A}}}1-\frac{1}{y_0}\sum_{\substack{x<n\leq x+y_0\\ n\in\mathcal{A}}}1\right\vert\mathrm{d}x\\
\leq\frac{1}{X}\int_X^{2X}\left\vert\frac{1}{h}\sum_{\substack{x<n\leq x+h\\ n\in\mathcal{A}\cap{\mathcal{S}}}}1-\frac{1}{y_0}\sum_{\substack{x<n\leq x+y_0\\ n\in\mathcal{A}\cap{\mathcal{S}}}}1\right\vert\mathrm{d}x+\frac{2}{X}\sum_{\substack{X<n\leq 2X+y_0 \\ n\in\mathcal{A}\smallsetminus{\mathcal{S}}}}1.
\end{multline*}
Grâce à la Proposition $\ref{estimationintegrale}$ et à l'hypothèse~\ref{hyp7}, en utilisant l'inégalité de Cauchy-Schwarz on obtient que l'intégrale de droite est un $o\left(\delta(X)\right)$. Par ailleurs, on a 
\[\frac{1}{X}\sum_{\substack{X<n\leq 2X+y_0 \\ n\in\mathcal{A}\smallsetminus{\mathcal{S}}}}1\leq\frac{1}{X}\sum_{\substack{n\sim X\\ n\in\mathcal{A}\smallsetminus{\mathcal{S}}}}1+\frac{y_0}{X}\leq\frac{1}{X}\left\vert\left(\mathcal{A}\smallsetminus\mathcal{S}\right)\cap\left[X,2X\right]\right\vert+o\left(\delta(X)\right)\]
avec l'hypothèse \ref{hyp1}. On déduit des l'hypothèses~\ref{hyp2} et~\ref{hyp8} que cette dernière quantité est $o\left(\delta(X)\right)$. Ainsi, 
\[\frac{1}{X}\int_X^{2X}\left\vert\frac{1}{h}\sum_{\substack{x<n\leq x+h\\ n\in\mathcal{A}}}1-\frac{1}{y_0}\sum_{\substack{x<n\leq x+y_0\\ n\in\mathcal{A}}}1\right\vert\mathrm{d}x=o\left(\delta(X)\right).\]
On en déduit aisément le Théorème $\ref{th1}$.
\end{proof}

Afin de démontrer la Proposition $\ref{estimationintegrale}$, on utilise la borne de Parseval pour réduire le problème à la majoration d'une intégrale faisant inetervenir le polyôme de Dirichlet associé à l'ensemble $\mathcal{A}\cap{\mathcal{S}}$.

\subsection{Utilisation de la borne de Parseval}\label{sssectionparseval}

Grâce au Lemme $\ref{Parseval}$, la Proposition $\ref{estimationintegrale}$ se déduit facilement de la proposition suivante. Avec les notations introduites précédemment, on pose
\[B(s):=\sum_{\substack{n\sim X\\n\in\mathcal{A}\cap{\mathcal{S}}}}\frac{1}{n^s}\hspace{2cm}(s\in\C).\]
\begin{prop}\label{majintegrale}
Sous les hypothèses~\ref{hyp1} à~\ref{hyp6} du Théorème~$\ref{th1}$, lorsque $T_0< T\leq X$ et $Q_1$ est suffisamment grand, on a
\begin{multline*}\int_{T_0}^T\left\vert B(1+it)\right\vert^2\mathrm{d}t\\
\ll\left(\frac{TQ_1}{\delta(X)X}+\theta(X)\right)\frac{\delta(X)^2}{H_1}+H_{\infty}^2(\log X)^{3-\varepsilon}\left(\frac{(\log X)^2}{T_0^2}+P_{\infty}^{-(\log X)^{-2/3-\varepsilon/3}}\right)
\end{multline*}
lorsque $X$ tend vers l'infini.
\end{prop}
En effet, la borne fournie par le Lemme~$\ref{MVT}$ est directement suffisante lorsque $T>X$. Il reste désormais à appliquer la méthode de découpage des polynômes de Dirichlet permettant de prouver cette proposition. 

\subsection{Factorisation du polynôme de Dirichlet associé à $\mathcal{A}\cap{\mathcal{S}}$}

La factorisation du polynôme $B$ s'effectue à l'aide du lemme suivant.

\begin{lemme}\label{decompo}
Avec les notations précédentes, lorsque $\mathcal{A}$ vérifie l'hypothèse~\ref{hyp4} du Théorème~$\ref{th1}$, et que $X$ est suffisamment grand, on a pour tout $j\in\left[1,J\right]\cup\{\infty\}$ et tout $s\in\C$,
\begin{equation}\label{facto}
B(s)=\sum_{\substack{n\sim X\\n\in\mathcal{A}\cap{\mathcal{S}}}}\frac{1}{n^s}=\sum_{v\in\mathcal{I}_j}Q_{v,H_j}(s)R_{v,H_j}(s)+N_{H_j}(s),
\end{equation}
où l'on a posé
\begin{align*}
Q_{v,H_j}(s)&:=\sum_{\substack{\e^{v/H_j}< p\leq\e^{(v+1)/H_j} \\ P_j< p\leq Q_j}}\frac{1}{p^s},\\
R_{v,H_j}(s)&:=\sum_{\substack{m\sim X\e^{-v/H_j} \\ m\in\mathcal{A}' \\ \left(m,\prod_{\e^{v/H_j}< p\leq \e^{(v+1)/H_j}}p\right)=1}}\frac{1}{m^s}\frac{1}{\omega_{]P_j,Q_j]}(m)+1},\\
N_{H_j}(s)&:=\sum_{\substack{2X<n\leq 2X\e^{1/H_j} \\ n\in\mathcal{A}\cap\mathcal{S}}}\frac{c_{n,j}}{n^s},
\end{align*}
pour certains complexes $c_{n,j}$ bornés en module par $1$.
\end{lemme}

\begin{proof}
Dans un premier temps, on montre qu'il existe des complexes $c_{n,j}$ tels que pour tout $s\in\C$,
\begin{equation}\label{eqcemlqskd}
\sum_{v\in\mathcal{I}_j}Q_{v,H_j}(s)R_{v,H_j}(s)=\sum_{n\in\mathcal{A}\cap{\mathcal{S}}}\frac{c_{n,j}}{n^s}.
\end{equation}
Soit $v\in\mathcal{I}_j$, puis $p$ et $m$ deux entiers apparaissant respectivement dans les sommes définissant $Q_{v,H_j}$ et $R_{v,H_j}$. Par définition de $\mathcal{A}'$, il existe $n\in\mathcal{A}\cap\mathcal{S}\cap[X;2X]$ et $q\in\bigcup_{1\leq i\leq J}]P_i,Q_i]$ tels que $m=\frac{n}{q}$. Alors $q\in\left[\frac{\e^{v/H_j}}{2},2\e^{v/H_j}\right]$ et donc, lorsque $X$ est assez grand, on obtient $q\in]P_j,Q_j]$. Par ailleurs, $mq=n\in\mathcal{A}\cap\mathcal{S}$ et $p\nmid m$ par définition de $R_{v,H_j}$. On a alors $mp\in\mathcal{A}$ avec l'hypothèse~\ref{hyp4}. De plus, $p\in\mathcal{S}_j$ est premier avec $m\in\cap_{\substack{i\in[1,J]\cup\{\infty\}\\ i\neq j}}\mathcal{S}_i$, donc $mp\in\mathcal{S}$, ce qui fournit~$(\ref{eqcemlqskd})$. Il est ensuite facile de voir que $X<mp\leq 2X\e^{1/H_j}$, et donc $c_{n,j}=0$ si $n\leq X$ ou $n>2X\e^{1/H_j}$. On conclut alors en faisant l'observation suivante. Pour tout $X<n\leq 2X\e^{1/H_j}$, il y a au plus $\omega_{]P_j,Q_j]}(n)$ décompositions possibles de $n$ sous la forme $mp$ précédente, et il y en a exactement $\omega_{]P_j,Q_j]}(n)$ lorsque $n\sim X$.
\end{proof}

\subsection{Majoration de l'intégrale}\label{majorationdelintegrale}

\begin{proof}[Démonstration de la Proposition $\ref{majintegrale}$]
Soit $T_0<T\leq X$. On met en place le découpage de l'intervalle $\left[T_0,T\right]$. Soit, pour $1\leq j\leq J$ (on rappelle que $\varepsilon$ intervient dans la définition de $\mathcal{S}$),
\begin{equation*}
\alpha_j:=\frac{1}{4}-\varepsilon\left(1+\frac{1}{2j}\right).
\end{equation*}
On a
\begin{equation*}
\frac{1}{4}-\frac{3}{2}\varepsilon=\alpha_1\leq \alpha_2\leq ...\leq \alpha_J\leq\frac{1}{4}-\varepsilon.
\end{equation*}
On écrit alors 
\begin{equation}\label{unionT0T}
\left[T_0,T\right]=\bigsqcup_{j=1}^J\mathcal{T}_j\sqcup\mathcal{U}
\end{equation}
où $t\in\mathcal{T}_j$ si $j$ est le plus petit indice tel que pour tout $v\in\mathcal{I}_j$ on ait
\begin{equation}\label{inegQ}
\left\vert Q_{v,H_j}(1+it)\right\vert\leq\left(\e^{v/H_j}\right)^{-\alpha_j},
\end{equation}
et $t\in\mathcal{U}$ si $t$ ne vérifie ces conditions pour aucun $1\leq j\leq J$. On traite désormais les intégrales sur les $\mathcal{T}_j$ et sur $\mathcal{U}$ de manière assez similaire à~\cite{matoradzi}.\\

\textbf{Cas $j=1$ : }\label{casj=1debut}Avec le Lemme~$\ref{decompo}$ et l'inégalité de Cauchy-Schwarz, on a
\[\int_{\mathcal{T}_1}\left\vert B(1+it)\right\vert^2\mathrm{d}t\ll\left\vert\mathcal{I}_1\right\vert\sum_{v\in\mathcal{I}_1}\int_{\mathcal{T}_1}\left\vert Q_{v,H_1}(1+it)R_{v,H_1}(1+it)\right\vert^2\mathrm{d}t+\int_{\mathcal{T}_1}\left\vert N_{H_1}(1+it)\right\vert^2\mathrm{d}t.\]
Avec le Lemme~$\ref{IMVT}$ et les hypothèses~\ref{hyp1}~et~\ref{hyp3}, on a
\begin{equation*}
\int_{\mathcal{T}_1}\left\vert N_{H_1}(1+it)\right\vert^2\mathrm{d}t\ll\left(\frac{T}{\delta(X)X}+\theta(X)\right)\frac{\delta(X)^2}{H_1}.
\end{equation*}
Par ailleurs, avec la définition de $\mathcal{T}_1$, le Lemme~$\ref{IMVT}$ et l'hypothèse~\ref{hyp5}, on a
\begin{align*}
\left\vert\mathcal{I}_1\right\vert\sum_{v\in\mathcal{I}_1}\int_{\mathcal{T}_1}\left\vert Q_{v,H_1}(1+it)R_{v,H_1}(1+it)\right\vert^2\mathrm{d}t&\ll\left\vert\mathcal{I}_1\right\vert\sum_{v\in\mathcal{I}_1}\e^{-2v\alpha_1/H_1}\left(\frac{T\e^{v/H_1}}{\delta(X)X}+\theta(X)\right)\delta(X)^2\\
&\ll\left(\frac{TQ_1}{\delta(X)X}+\theta(X)\right)\delta(X)^2H_1^2(\log Q_1)P_1^{-2\alpha_1}\\
&\ll\left(\frac{TQ_1}{\delta(X)X}+\theta(X)\right)\frac{\delta(X)^2}{H_1}.
\end{align*}
\\
\indent\textbf{Cas $2\leq j\leq J$ : }\label{methodeth1}Pour tout $t\in\mathcal{T}_j$, il existe au moins un $u\in\mathcal{I}_{j-1}$ tel que~$(\ref{inegQ})$ ne soit pas vérifiée, et donc $\left\vert Q_{u,H_{j-1}}(1+it)\left(\e^{u/H_{j-1}}\right)^{2\alpha_{j-1}}\right\vert >1$. On décompose donc $\mathcal{T}_j$ selon ces indices $u$, de la manière suivante,
\[\mathcal{T}_j=\bigcup_{u\in\mathcal{I}_{j-1}}\mathcal{T}_{j,u}.\]
De la même manière que précédemment, on a
\[\int_{\mathcal{T}_j}\left\vert B(1+it)\right\vert^2\mathrm{d}t\ll\left\vert\mathcal{I}_j\right\vert\sum_{v\in\mathcal{I}_j}\int_{\mathcal{T}_j}\left\vert Q_{v,H_j}(1+it)R_{v,H_j}(1+it)\right\vert^2\mathrm{d}t+\int_{\mathcal{T}_j}\left\vert N_{H_j}(1+it)\right\vert^2\mathrm{d}t\]
où la dernière intégrale est $\ll\left(\frac{T}{\delta(X)X}+\theta(X)\right)\frac{\delta(X)^2}{H_j}$. En prenant le terme maximal dans la somme, il existe $u_{j-1}\in\mathcal{I}_{j-1}$ et $v_j\in\mathcal{I}_j$ tels que
\[\left\vert\mathcal{I}_j\right\vert\sum_{v\in\mathcal{I}_j}\int_{\mathcal{T}_j}\left\vert Q_{v,H_j}(1+it)R_{v,H_j}(1+it)\right\vert^2\mathrm{d}t\ll\left\vert\mathcal{I}_j\right\vert^2\left\vert\mathcal{I}_{j-1}\right\vert\e^{-2v_j\alpha_j/H_j}\int_{\mathcal{T}_{j,u_{j-1}}}\left\vert R_{v_j,H_j}(1+it)\right\vert^2\mathrm{d}t.\]
Si l'on appliquait directement le Lemme~$\ref{IMVT}$ à la dernière intégrale, le polynôme $R_{v_j,H_j}$ étant trop court, on perdrait un facteur potentiellement de la taille de $Q_j$. Cette perte serait trop importante, à moins que $T$ soit assez petit, ce qui reviendrait à s'intéresser à de plus grands intervalles dans le problème initial. Pour contourner ce problème, Matomäki et Radziwi\l\l\ pensent à utiliser la définition de $\mathcal{T}_{j,u}$ afin d'agrandir artificiellement la taille du polynôme dans l'intégrale. Plus précisément, pour tout $\ell\geq 1$ on a
\[\int_{\mathcal{T}_{j,u_{j-1}}}\left\vert R_{v_j,H_j}(1+it)\right\vert^2\mathrm{d}t\leq\e^{2\ell u_{j-1}\alpha_{j-1}/H_{j-1}}\int_{\mathcal{T}_{j,u_{j-1}}}\left\vert Q_{u_{j-1},H_{j-1}}(1+it)^{\ell}R_{v_j,H_j}(1+it)\right\vert^2\mathrm{d}t.\]
On pose $Y_1:=\e^{u_{j-1}/H_{j-1}}$ et $Y_2:=\e^{v_j/H_j}$. On choisit de prendre $\ell:=\left\lceil\frac{\log Y_2}{\log Y_1}\right\rceil=\left\lceil\frac{v_j/H_j}{u_{j-1}/H_{j-1}}\right\rceil$, de sorte que le polynôme résultant soit supporté dans un intervalle de taille assez proche de $X$. Pour être plus précis, pour tout $s\in\C$ on a
\[Q_{u_{j-1},H_{j-1}}(s)^{\ell}R_{v_j,H_j}(s)=\sum_{X\leq n\leq 2^{\ell +1}Y_1X}\frac{a_n}{n^s}\]
pour certains complexes vérifiant $\left\vert a_n\right\vert\leq\left(\ell+1\right)!$, puisque tout $m\in\mathcal{A}'$ a au plus $1$ seul facteur premier dans l'intervalle $\left]\e^{u_{j-1}/H_{j-1}},\e^{(u_{j-1}+1)/H_{j-1}}\right]$. Par ailleurs, pour que $a_n$ soit non nul, il faut au moins que $n$ soit de la forme $p_1\cdots p_{\ell}m$ avec $p_i\sim Y_1$ pour tout $1\leq i\leq\ell$, et $m\in\mathcal{A}'$ vérifie $m\sim\frac{X}{Y_2}$. Avec le Lemme~$\ref{IMVT}$, on a donc
\begin{multline*}
\int_{\mathcal{T}_{j,u_{j-1}}}\left\vert Q_{u_{j-1},H_{j-1}}(1+it)^{\ell}R_{v_j,H_j}(1+it)\right\vert^2\mathrm{d}t\ll T\left(\ell+1\right)!^2\sum_{p_1, ..., p_{\ell}\sim Y_1}\frac{1}{(p_1\cdots p_{\ell})^2}\sum_{\substack{m\sim\frac{X}{Y_2} \\ m\in\mathcal{A}'}}\frac{1}{m^2}\\
+T\left(\ell+1\right)!^2\sum_{1\leq h\leq\frac{2^{\ell +1}Y_1X}{T}}\sum_{p_1, ..., p_{\ell}\sim Y_1}\frac{1}{(p_1\cdots p_{\ell})^2}\sum_{\substack{m\sim\frac{X}{Y_2} \\ m\in\mathcal{A}' \\ p_1\cdots p_{\ell}m+h\in\bigcup_{q_1, ..., q_{\ell}\sim Y_1}q_1\cdots q_{\ell}\mathcal{A}'}}\frac{1}{m^2}.
\end{multline*}
En utilisant les hypothèses~\ref{hyp5}~et~\ref{hyp6}, on obtient alors
\[\int_{\mathcal{T}_{j,u_{j-1}}}\left\vert Q_{u_{j-1},H_{j-1}}(1+it)^{\ell}R_{v_j,H_j}(1+it)\right\vert^2\mathrm{d}t\ll\left(\frac{T}{\delta(X)X}+\theta(X)\right)\delta(X)^2\left(\ell !\right)^{O(1)}Y_1.\]
On a donc
\begin{multline}\label{huid}
\left\vert\mathcal{I}_j\right\vert\sum_{v\in\mathcal{I}_j}\int_{\mathcal{T}_j}\left\vert Q_{v,H_j}(1+it)R_{v,H_j}(1+it)\right\vert^2\mathrm{d}t\\
\ll\left(\frac{T}{\delta(X)X}+\theta(X)\right)\delta(X)^2\left(\ell !\right)^{O(1)}\left\vert\mathcal{I}_j\right\vert^2\left\vert\mathcal{I}_{j-1}\right\vert Q_{j-1}^{1+2\alpha_{j-1}}\e^{-2(\alpha_j-\alpha_{j-1})v_j/H_j}.
\end{multline}
Par ailleurs, on a
\[\ell\log\ell\leq\frac{v_j}{H_j}\frac{\log_2 Q_j}{\log P_{j-1}-1}+\log_2 Q_j+1.\]
Ainsi, comme $\alpha_j-\alpha_{j-1}\geq\frac{\varepsilon}{2j^2}$, lorsque $Q_1$ est suffisamment grand, uniformément pour $2\leq j\leq J$ le produit des termes après $\delta(X)^2$ dans~$(\ref{huid})$ est
\begin{align*}
&\ll (H_j\log Q_j)^2H_{j-1}(\log Q_{j-1})Q_{j-1}^{1+2\alpha_{j-1}}\exp\left(-\frac{v_j}{H_j}\left(\frac{\varepsilon}{j^2}-O\left(\frac{\log_2 Q_j}{\log P_{j-1}-1}\right)\right)\right)(\log Q_j)^{O(1)}\\
&\ll j^6P_1^{1/2}Q_{j-1}^{O(1)}\exp\left(-\frac{v_j}{H_j}\left(\frac{\varepsilon}{j^2}-O\left(\frac{\log_2 Q_j}{\log P_{j-1}-1}\right)\right)\right)\\
&\ll j^6P_1^{1/2}Q_{j-1}^{O(1)}P_j^{-\varepsilon/(2j^2)}\ll \frac{1}{j^2P_1}\ll\frac{1}{j^2H_1}.
\end{align*}
Finalement, on a 
\begin{equation*}
\sum_{j=2}^J\int_{\mathcal{T}_j}\left\vert B(1+it)\right\vert^2\mathrm{d}t\ll\left(\frac{T}{\delta(X)X}+\theta(X)\right)\frac{\delta(X)^2}{H_1}.
\end{equation*}
Comme pour le cas $j=1$, l'hypothèse~\ref{hyp7} donne alors une bonne majoration.
\\
\indent\textbf{Cas de $\mathcal{U}$ : }Cette fois, on écrit
\[\int_{\mathcal{U}}\left\vert B(1+it)\right\vert^2\mathrm{d}t\ll\left\vert\mathcal{I}_{\infty}\right\vert\sum_{v\in\mathcal{I}_{\infty}}\int_{\mathcal{U}}\left\vert Q_{v,H_{\infty}}(1+it)R_{v,H_{\infty}}(1+it)\right\vert^2\mathrm{d}t+\int_{\mathcal{U}}\left\vert N_{H_{\infty}}(1+it)\right\vert^2\mathrm{d}t.\]
De manière analogue aux points précédents,
\[\int_{\mathcal{U}}\left\vert N_{H_{\infty}}(1+it)\right\vert^2\mathrm{d}t\ll\left(\frac{T}{\delta(X)X}+\theta(X)\right)\frac{\delta(X)^2}{H_{\infty}}.\]
Soit $u\in\mathcal{I}_{\infty}$ un indice maximisant l'intégrale dans la somme précédente. Soit $\mathcal{T}\subset\mathcal{U}$ un ensemble $1$-espacé tel que l'on ait
\[\int_{\mathcal{U}}\left\vert Q_{u,H_{\infty}}(1+it)R_{u,H_{\infty}}(1+it)\right\vert^2\mathrm{d}t\ll\sum_{t\in\mathcal{T}}\left\vert Q_{u,H_{\infty}}(1+it)R_{u,H_{\infty}}(1+it)\right\vert^2.\]
Avec le Lemme~$\ref{majlongpol}$, pour tout $t\in\mathcal{T}$, on a
\[Q_{u,H_{\infty}}(1+it)\ll\frac{\log X}{T_0}+P_{\infty}^{-(\log X)^{-2/3-\varepsilon/3}}.\]
On a alors
\begin{multline*}
\left\vert\mathcal{I}_{\infty}\right\vert\sum_{v\in\mathcal{I}_{\infty}}\int_{\mathcal{U}}\left\vert Q_{v,H_{\infty}}(1+it)R_{v,H_{\infty}}(1+it)\right\vert^2\mathrm{d}t\\
\ll H_{\infty}^2\left(\log Q_{\infty}\right)^2\left(\frac{(\log X)^2}{T_0^2}+P_{\infty}^{-(\log X)^{-2/3-\varepsilon/3}}\right)\sum_{t\in\mathcal{T}}\left\vert R_{u,H_{\infty}}(1+it)\right\vert^2
\end{multline*}
En majorant trivialement la dernière somme avec le Lemme~$\ref{halaszentiers}$, on trouve
\[\sum_{t\in\mathcal{T}}\left\vert R_{u,H_{\infty}}(1+it)\right\vert^2\ll\left(\frac{\left\vert\mathcal{T}\right\vert\sqrt{T}Q_{\infty}}{X}+1\right)(\log X)\ll\log X,\]
où la dernière inégalité est obtenue grâce au Lemme~$\ref{grandesvaleurspolynome}$, en se rappelant la définition de $\mathcal{U}$ et le fait que $T\leq X$. On atteint ainsi la conclusion désirée.
\end{proof}

\section{Application aux entiers ayant $k\asymp\log_2 X$ facteurs premiers}\label{applis}

Dans cette section, on démontre le Théorème~\ref{corollaire}. Pour cela, on commence par énoncer quelques propriétés classiques des entiers ayant exactement $k$ facteurs premiers. On démontre ensuite un théorème d'indépendance entre $\omega(a_1n+b_1)$ et $\omega(a_2n+b_2)$ à la sous-section~\ref{ssseccribidim}. On applique enfin le Théorème~\ref{th1}.

Pour $x\geq 3$ et $k\geq 1$, on rappelle les notations
\begin{equation*}\label{formuleexpldelt}
\delta_k(x):=\lambda\left(\kappa\right)\frac{(\log_2 x)^{k-1}}{(\log x)(k-1)!},\hspace{2cm}\kappa:=\frac{k-1}{\log_2 x},
\end{equation*}
où
\begin{equation}\label{deffctla}
\lambda(z):=\frac{1}{\Gamma(z+1)}\prod_{p\geq 2}\left(1+\frac{z}{p-1}\right)\left(1-\frac{1}{p}\right)^z.
\end{equation}
En particulier, lorsque $0<r<R$ sont fixés et $r<\kappa\leq R$, on a
\begin{equation}\label{tyuiop}
\delta_k(x)\asymp\frac{1}{(\log x)^{Q(\kappa)}\sqrt{\log_2 x}},
\end{equation}
où $Q(\kappa):=\kappa\log\kappa-\kappa+1$.

\subsection{Propriétés classiques}

\begin{lemme}\label{hyp1application}
Il existe une constante $c>0$ vérifiant l'énoncé suivant. Soient $A, R, \beta_1>0$ et $0<\beta_2<2$ fixés. Soit, pour $x\geq 3$, $y=y(x)\geq 3$ tel que $\log_2 y\leq A\sqrt{\log_2 x}$. Uniformément pour $x\geq 3$, $x^{1-\frac{c}{\log y}}<x'\leq x$, $1\leq k\leq R\log_2 x$, et $f$ une fonction multiplicative telle que
\[0\leq f(p^{\nu})\leq\beta_1\beta_2^{\nu-1}\hspace{2cm}(p\geq 2,\ \nu\geq 1),\]
et $f(p^{\nu})=1$ dès que $p>y$, on a
\begin{equation}\label{formulehypap1}
\sum_{\substack{x<n\leq x+x' \\ n\in\mathcal{E}_k}}f(n)=\delta_k(x)x'\prod_{p\leq x}\frac{1+\kappa\sum_{\nu=1}^{\infty}\frac{f(p^{\nu})}{p^{\nu}}}{1+\frac{\kappa}{p-1}}\left\{1+O\left(\frac{k(\log_2 y)^2}{(\log_2 x)^2}\right)\right\},
\end{equation}
lorsque $x$ tend vers l'infini. En particulier, il existe une constante absolue $\delta>0$ telle que pour tous $x^{1-\delta}<x'\leq x$ et $1\leq k\leq R\log_2 x$ on ait
\begin{equation*}
\pi_k(x+x')-\pi_k(x)=\delta_k(x)x'\left(1+o(1)\right),
\end{equation*}
lorsque $x$ tend vers l'infini.
\end{lemme}

Dans le cas $x'=x$, ce lemme est énoncé sans démonstration par Tenenbaum~\cite{tenenbaumbillingsley}, comme une conséquence facile de la méthode de Selberg-Delange. On poursuit ici cette idée, en utilisant les travaux de Cui et Wu~\cite{cuiwu}, qui ont adapté la méthode de Selberg-Delange aux petits intervalles.
\begin{proof}

On suppose donnés les paramètres de l'énoncé. Sans perte de généralité, on suppose $R,\beta_2 >1$. On note que pour tout $\sigma>\frac{\log\beta_2}{\log 2}$, on a
\[\sum_{\nu\geq 1}\frac{f(p^{\nu})}{p^{\sigma\nu}}\leq\frac{\beta_1}{p^{\sigma}-\beta_2}\ll\frac{1}{p^{\sigma}}.\]
Pour $\vert z\vert\leq R$ et $0\leq t'\leq t$, on pose
\[S(t,t'):=\sum_{t< n\leq t+t'}f(n)z^{\omega(n)}.\]
On estime $S(x,x')$ pour ensuite démontrer le lemme via la formule de Cauchy. Pour $z\in\C$, on définit la fonction multiplicative $\tau_z$ sur les puissances de nombres premiers par
\[\tau_z(p^{\nu}):=\binom{z+\nu-1}{\nu},\]
de sorte que pour tout $s\in\C$ de partie réelle $\sigma>1$ on ait
\[\zeta(s)^z=\prod_{p\geq 2}\left(1-\frac{1}{p^s}\right)^{-z}=\sum_{n\geq 1}\frac{\tau_z(n)}{n^s}.\]
On définit alors implicitement la fonction multiplicative $g_z$ par l'égalité de convolution\begin{equation}\label{efconvolkui}
fz^{\omega(\cdot)}=\tau_z\ast g_z.
\end{equation}
Pour $\sigma>\max\left(\frac{\log\beta_2}{\log 2},\frac{1}{2}\right)$, on a donc
\[\sum_{n\geq 1}\frac{g_z(n)}{n^s}=\prod_{p\geq 2}\left(1+z\sum_{\nu\geq 1}\frac{f(p^{\nu})}{p^{s\nu}}\right)\left(1-\frac{1}{p^s}\right)^{z},\]
et
\[g_z=fz^{\omega(\cdot)}\ast\tau_{-z}.\]
En particulier, $g_z(p)=0$ lorsque $p>y$. Soit $\alpha:=\frac{c_1}{\log y}$ où $c_1>0$ est une constante suffisamment petite pour que $2^{1-\alpha}>\beta_2$ pour tout $x$. On remarque tout de suite que $\log y=o\left((\log x)^{o(1)}\right)$, et que pour tous $p\geq 2,\nu\geq 1$, on a
\[\vert g_z(p^{\nu})\vert\ll\beta_2^{\nu}.\]
Ainsi, uniformément pour tout $t\geq 1$, on a
\begin{align}
\sum_{d>t}\frac{\vert g_z(d)\vert}{d}&\leq t^{-\alpha}\sum_{P^+(d)\leq y}\frac{\vert g_z(d)\vert}{d^{1-\alpha}}\prod_{p>y}\left(1+\sum_{\nu\geq 2}\frac{\vert g_z(p^{\nu})\vert}{p^{\nu}}\right)\nonumber\\
&\ll t^{-\alpha}\prod_{2\leq p\leq y}\left(1+\sum_{\nu\geq 1}\frac{\vert g_z(p^{\nu})\vert}{p^{(1-\alpha)\nu}}\right)\prod_{p>y}\left(1+\sum_{\nu\geq 2}\frac{\vert g_z(p^{\nu})\vert}{p^{\nu}}\right)\nonumber\\
&\ll\exp\left(-c_1\frac{\log t}{\log y}\right)(\log x)^{o(1)}.\label{hjfngjhnjgbhgf}
\end{align}
Soient $c_2<1$ une constante suffisamment petite et $w:=x^{c_2}$. En définissant
\begin{equation*}
S_1(z):=\sum_{\substack{d\leq w}}g_z(d)\sum_{\frac{x}{d}<m\leq\frac{x+x'}{d}}\tau_z(m),\hspace{1cm}S_2(z):=\sum_{\substack{w<d\leq 2x}}g_z(d)\sum_{\substack{\frac{x}{d}<m\leq\frac{x+x'}{d}}}\tau_z(m),
\end{equation*}
on a
\begin{align*}
S(x,x')&=\sum_{d\leq 2x}g_z(d)\sum_{\frac{x}{d}<m\leq\frac{x+x'}{d}}\tau_z(m)\\
&=S_1(z)+S_2(z).
\end{align*}
On majore $\vert S_2(z)\vert$ avec la méthode de Rankin. On a
\begin{align}
\vert S_2(z)\vert&\leq\sum_{\substack{w<d\leq 2x}}\vert g_z(d)\vert\sum_{\frac{x}{d}<m\leq\frac{x+x'}{d}}\vert \tau_z(m)\vert\frac{2x}{dm}\nonumber\\
&\ll x(\log x)^R\sum_{\substack{d>w}}\frac{\vert g_z(d)\vert}{d}\nonumber\\
&\ll x(\log x)^{R+o(1)}\exp\left(-c_1c_2\frac{\log x}{\log y}\right).\label{eswxszsdfgbhj}
\end{align}
La dernière inégalité est obtenue en utilisant~$(\ref{hjfngjhnjgbhgf})$. Lorsque les constantes $c_1, c_2$ puis $c$ sont suffisamment petites, d'après~\cite[theorem~1.1]{cuiwu}, uniformément pour $1\leq d\leq w$, on a
\[\sum_{\frac{x}{d}<m\leq\frac{x+x'}{d}}\tau_z(m)=\frac{x'}{d}\left(\log\frac{x}{d}\right)^{z-1}\left(\frac{1}{\Gamma(z)}+O\left(\frac{1}{\log x}\right)\right).\]
Il s'ensuit
\[S_1(z)=x'(\log x)^{z-1}\sum_{\substack{d\leq w}}\frac{g_z(d)}{d}\left(1-\frac{\log d}{\log x}\right)^{z-1}\left(\frac{1}{\Gamma(z)}+O\left(\frac{1}{\log x}\right)\right).\]
Soit $w':=y^{\frac{\log_2 x}{c_1}}$. On traite la somme ci-dessus en séparant les $d\leq w'$, pour lesquels $\left(1-\frac{\log d}{\log x}\right)^{z-1}=1+O\left((\log x)^{-1+o(1)}\right)$, des entiers $d$ tels que $w'<d\leq w$. Ce faisant et en utilisant~$(\ref{hjfngjhnjgbhgf})$, on obtient
\[S_1(s)=x'(\log x)^{z-1}\left(\frac{1}{\Gamma(z)}\sum_{d\leq w'}\frac{g_z(d)}{d}+O\left(\frac{1}{(\log x)^{1-o(1)}}\right)\right).\]
Finalement, en utilisant de nouveau~$(\ref{hjfngjhnjgbhgf})$ pour éliminer la condition $d\leq w'$, puis~$(\ref{eswxszsdfgbhj})$ et l'hypothèse sur $x'$, on a
\begin{equation}\label{sommeSevaluee}
S(x,x')=x'(\log x)^{z-1}\left\{\lambda_0(z)+O\left(\frac{1}{(\log x)^{1-o(1)}}\right)\right\},
\end{equation}
où l'on a posé $\lambda_0(z):=zh_0(z)$ et
\[h_0(z):=\frac{1}{\Gamma(z+1)}\prod_{p\geq 2}\left(1+z\sum_{\nu\geq 1}\frac{f(p^{\nu})}{p^{\nu}}\right)\left(1-\frac{1}{p^s}\right)^z.\]

La fonction $h_0$ est entière. Pour $p\geq 2$, on pose
\[s_p:=(p-1)\sum_{\nu\geq 1}\frac{f(p^{\nu})}{p^{\nu}},\]
de sorte que $0\leq s_p\leq\beta_1\frac{p-1}{p-\beta_2}$ pour $p\geq 2$, et $s_p=1$ pour $p>y$. Un calcul élémentaire fournit d'une part
\begin{equation}\label{egalhsp}
h_0(z)=\frac{1}{\Gamma(z+1)}\prod_{p\geq 2}\left(1+\frac{s_pz}{p-1}\right)\left(1-\frac{1}{p}\right)^z,
\end{equation}
et d'autre part
\begin{multline}\label{hseconde}
h_0''(z)=h_0(z)\left\{\left(\frac{-\Gamma'}{\Gamma}\right)'\!\!(z+1)-\sum_{p\geq 2}\frac{s_p^2}{(p-1+s_pz)^2}\right.\\
\left.+\left(\frac{-\Gamma'}{\Gamma}(z+1)+\sum_{p\geq 2}\frac{s_p}{p-1+s_pz}+\log\left(1-\frac{1}{p}\right)\right)^2\right\},
\end{multline}
où tous les pôles du terme entre accolades sont compensés par les zéros de $h_0$. Avec~$(\ref{egalhsp})$, on voit que $h_0$ ne s'annule pas dans le demi-plan $\sigma>-\frac{2-\beta_2}{\beta_1}$. De~$(\ref{egalhsp})$ et~$(\ref{hseconde})$, on déduit qu'uniformément pour tous $\vert z\vert, \vert z'\vert\leq R$ de partie réelle $\sigma, \sigma'\geq 0$, on a
\begin{align}
\left\vert\frac{h_0(z)}{h_0(z')}\right\vert&\ll(\log y)^{O(\vert\sigma-\sigma'\vert)},\label{inegh4}\\
\vert h_0''(z)\vert&\ll(\log_2 y)^2\vert h_0(z)\vert,\label{inegh3}\\
\vert h_0''(-z)\vert&\ll(\log_2 y)^2(\log y)^{O(\sigma)}.\label{ineg5h}
\end{align}

D'après la formule de Cauchy, le terme de gauche de~$(\ref{formulehypap1})$ vaut, pour tout $r>0$,
\begin{equation}\label{bhghgbhb}
\frac{1}{2i\pi}\int_{\vert z\vert =r}\frac{S(x,x')}{z^{k+1}}\mathrm{d}z.
\end{equation}
On estime cette intégrale avec l'expression~$(\ref{sommeSevaluee})$ de $S(x,x')$, qui dépend bien sûr de la variable~$z$. On commence par le terme principal. Dans le cas $k=1$, l'intégrale du terme principal de $S(x,x')$ vaut $\frac{x'}{\log x}h_0(0)$, comme voulu. Dans le cas $k\geq 2$, on pose $r:=\frac{k-1}{\log_2 x}=\kappa$, et on estime
\begin{equation}\label{ksjdhfdksjh}
\int_{\vert z\vert =r}\frac{(\log x)^zh_0(z)}{z^k}\mathrm{d}z=h_0(r)\int_{\vert z\vert =r}\frac{(\log x)^z}{z^k}\mathrm{d}z+\int_{\vert z\vert =r}\frac{(\log x)^z(h_0(z)-h_0(r))}{z^k}\mathrm{d}z.
\end{equation}
La première intégrale du membre de droite vaut $2i\pi\frac{(\log_2 x)^{k-1}}{(k-1)!}$, ce qui fournit le terme principal de~$(\ref{formulehypap1})$. En effet, avec~$(\ref{egalhsp})$, on a
\[h_0(r)=\lambda(r)\prod_{p\leq x}\frac{1+r\sum_{\nu=1}^{\infty}\frac{f(p^{\nu})}{p^{\nu}}}{1+\frac{r}{p-1}},\]
où la fonction $\lambda$ est définie par~$(\ref{deffctla})$. Comme dans~\cite[\chrom{2}.6.1]{tenenbaum}, on remarque par un calcul explicite que
\[\int_{\vert z\vert=r}(z-r)\frac{(\log x)^z}{z^k}\mathrm{d}z=0.\]
Ainsi, la dernière intégrale de~$(\ref{ksjdhfdksjh})$ vaut
\begin{equation*}\label{fcgfvbhnk}
\int_{\vert z\vert =r}\frac{(\log x)^z}{z^k}\left\{h_0(z)-h_0(r)-(z-r)h_0'(r)\right\}\mathrm{d}z.
\end{equation*}
Après utilisation de la formule de Taylor avec reste intégral, on majore
\[\int_{0}^1(1-t)\int_{\vert z\vert=r}\frac{(\log x)^z}{z^k}(z-r)^2h_0''(r+t(z-r))\mathrm{d}z\mathrm{d}t.\]
Pour $z$ de partie réelle $\sigma\leq 0$, on majore à l'aide de~$(\ref{inegh4})$ avec $z'=0$,~$(\ref{inegh3})$ et~$(\ref{ineg5h})$, obtenant ainsi
\begin{align*}
\int_{\substack{\vert z\vert=r \\ \sigma\leq 0}}\frac{(\log x)^z}{z^k}(z-r)^2h_0''(r+t(z-r))\mathrm{d}z&\ll r^{3-k}(\log_2 y)^2(\log y)^{O(r)}\\
&\ll\frac{(\log_2 x)^{k-1}}{(k-1)!}\frac{\e^{-k}k^{5/2}(\log_2 y)^2\e^{O\left(k/\sqrt{\log_2 x}\right)}}{(\log_2 x)^2}\\
&\ll\frac{(\log_2 x)^{k-1}}{(k-1)!}h_0(r)\frac{(\log_2 y)^2}{(\log_2 x)^2},
\end{align*}
où la dernière inégalité est obtenue en utilisant~$(\ref{inegh4})$ avec $z=0$. Lorsque $\sigma>0$, on majore grâce à~$(\ref{inegh3})$ puis~$(\ref{inegh4})$, obtenant ainsi
\begin{align*}
\int_{\substack{\vert z\vert=r \\ \sigma>0}}\!\!\frac{(\log x)^z}{z^k}(z-r)^2h_0''(r+t(z&-r))\mathrm{d}z\\
&\ll h_0(r)(\log_2 y)^2\int_{\substack{\vert z\vert=r \\ \sigma>0}}\!\!\frac{(\log x)^z}{z^k}(z-r)^2(\log y)^{O(\re(r-z))}\mathrm{d}z\\
&\ll\frac{h_0(r)(\log_2 y)^2}{r^{k-3}}\int_{-\frac{\pi}{2}}^{\frac{\pi}{2}}\e^{(k-1)\cos\theta}\left\vert\e^{i\theta}-1\right\vert^2\e^{O(r(1-\cos\theta)\log_2 y)}\mathrm{d}\theta.
\end{align*}
La dernière intégrale, par changement de variable $\cos\theta=u$, est
\begin{align*}
&\ll\e^{k-1}\int_0^1\e^{-(1-u)(k-1)\left(1+O(\log_2 y/\log_2 x)\right)}\sqrt{1-u}\mathrm{d}u\\
&\ll\e^{k-1}(k-1)^{3/2},
\end{align*}
la dernière inégalité étant obtenue en effectuant un second changement de variable, on posant $(1-u)(k-1)=v$. On obtient alors, uniformément pour $0\leq t\leq 1$,
\begin{equation*}
\int_{\substack{\vert z\vert=r}}\!\!\frac{(\log x)^z}{z^k}(z-r)^2h_0''(r+t(z-r))\mathrm{d}z\ll \frac{(\log_2 x)^{k-1}}{(k-1)!}h_0(r)\frac{k(\log_2 y)^2}{(\log_2 x)^2},
\end{equation*}
qui est convenable. En revenant à~$(\ref{bhghgbhb})$ puis~$(\ref{sommeSevaluee})$, il nous reste à majorer
\[\int_{\vert z\vert=r}\frac{(\log x)^{\re z}}{\vert z\vert^{k+1}}\mathrm{d}z,\]
ce qui correspond à~$(6.14)$ de~\cite[\chrom{2}.6.1]{tenenbaum}. On a alors
\[\int_{\vert z\vert=r}\frac{(\log x)^{\re z}}{\vert z\vert^{k+1}}\mathrm{d}z\ll\frac{(\log_2 x)^k}{k!},\]
qui est convenable.

\end{proof}

\begin{lemme}\label{hyp2application}
Soient $R>0$ fixé. Uniformément pour $x\geq 3$, $1\leq k\leq R\log_2 x$ et $2\leq P<Q\leq x$, on a
\begin{equation*}
\left\vert\left\{n\sim x\,:\quad\omega(n)=k,\ \left(n,\prod_{P<p\leq Q}p\right)=1\right\}\right\vert\ll\delta_k(x)x\left(\frac{\log P}{\log Q}\right)^{\kappa}.
\end{equation*}
De plus, lorsque $k=o\left(\sqrt{\log_2 x}\right)$, la constante implicite peut être remplacée par $1+o(1)$.
\end{lemme}
\begin{proof}
La première partie de l'énoncé est une conséquence directe de~\cite[lemme~1]{tenenbaumbillingsley}. On traite maintenant le cas où $k=o\left(\sqrt{\log_2 x}\right)$. Lorsque $k=1$, on a $\kappa=0$, et donc le théorème des nombres premiers fournit le résultat. Soit désormais $k\geq 2$ tel que $k=o\left(\sqrt{\log_2 x}\right)$. Pour tout entier $n\in\mathcal{E}_k\cap]x,2x]$, il existe $p^{\nu}\Vert n$ tel que $p^{\nu}> x^{1/k}$. Le cardinal qui nous intéresse est alors
\[\leq\frac{1}{(k-1)!}\sum_{\substack{p_1^{\nu_1}, ..., p_{k-1}^{\nu_{k-1}} \\ p_1^{\nu_1}\cdots p_{k-1}^{\nu_{k-1}}\leq 2x^{1-1/k} \\ p_i\notin]P,Q], \forall i\in [1,k-1]}}\sum_{p_k^{\nu_k}\sim\frac{x}{p_1^{\nu_1}\cdots p_{k-1}^{\nu_{k-1}}}}1.\]
Pour tous $0<\varepsilon<1$ et $0\leq t\leq 1-\varepsilon$, on a $\frac{1}{1-t}\leq 1+\frac{t}{\varepsilon}$. Ainsi, avec le théorème des nombres premiers sous la forme $\pi(x)=\frac{x}{\log x}+O\left(\frac{x}{(\log x)^2}\right)$, lorsque $x$ est suffisamment grand, il existe une constante $C>0$ telle que l'on ait
\[\sum_{p_k^{\nu_k}\sim\frac{x}{p_1^{\nu_1}\cdots p_{k-1}^{\nu_{k-1}}}}1\leq\frac{x}{p_1^{\nu_1}\cdots p_{k-1}^{\nu_{k-1}}}\left(\frac{1}{\log x}+2k\frac{\log\left(p_1^{\nu_1}\cdots p_{k-1}^{\nu_{k-1}}\right)}{(\log x)^2}+\frac{Ck^2}{(\log x)^2}\right).\]
On conclut alors aisément avec les formules de Mertens $\sum_{p^{\nu}\leq x}\frac{\log p^{\nu}}{p^{\nu}}=\log x+O(1)$ et $\sum_{p^{\nu}\leq x}\frac{1}{p^{\nu}}=\log_2 x+O(1)$, sachant que pour un tel $k$, on a $\delta_k(x)=(1+o(1))\frac{(\log_2 x)^{k-1}}{(\log x)(k-1)!}$.
\end{proof}

\subsection{Crible bi-dimensionnel pour les entiers ayant $k$ facteurs premiers}\label{ssseccribidim}

Étant donnés $\P$ un ensemble de nombres premier, $n\geq 1$ un entier, et $x, y\geq 1$ deux réels, on note
\begin{equation}\label{dEfEhetnP}
E(x):=1+\sum_{\substack{p\leq x\\ p\in\P}}\frac{1}{p},\hspace{2cm}n_y:=\prod_{\substack{p^{\nu}\Vert n \\ p\leq y}}p^{\nu}.
\end{equation}
On définit $\N_{\P}:=\left\{n\geq 1\,:\quad p\vert n\Rightarrow p\in\P\right\}$, l'ensemble des entiers ayant tous leurs facteurs premiers dans $\P$. Dans cette sous-section, on démontre le théorème suivant.

\begin{theorem}\label{nouvelleprop}
Soient $0<\varepsilon\leq\frac{1}{3}$ et $A\geq 2$ fixés. Il existe une constante $K=K(A,\varepsilon)>0$ vérifiant l'énoncé suivant. Soit $\P$ un ensemble de nombres premiers. Soient $x\geq 2$ un réel et $1\leq k_1,k_2\leq AE(x)$ des entiers. Soient, pour $i\in\{1,2\}$, $Q_i(X)=a_iX+b_i$ où $1\leq a_i\leq x^A$, $0\leq b_i\leq x^A$ et $(a_i,b_i)=1$, tels que $\Del\neq 0$. Lorsque $x^{\varepsilon}< x'\leq x$ on a
\begin{enumerate}[parsep=0cm, topsep=0.2cm]
\item\label{premiercascri} $\d\left\vert\left\{n\in]x,x+x']\,:\quad\omega(Q_1(n))=k_1,\ Q_1(n)\in\N_{\P}\right\}\right\vert\leq K\frac{a_1}{\varphi(a_1)}\frac{x'}{\log x}\frac{E(x)^{k_1-1}}{(k_1-1)!}$,
\item\label{deuxiemecascri} \mathtoolsset{firstline-afterskip=125pt}$\!\!\begin{multlined}[t]
\left\vert\left\{n\in]x,x+x']\,:\quad\forall i\in\{1,2\},\ \omega(Q_i(n))=k_i,\ Q_i(n)\in\N_{\P}\right\}\right\vert\\
\leq K\frac{\vert\Del\vert^{K} a_1a_2}{\varphi(\vert\Del\vert)^{K}\varphi(a_1)\varphi(a_2)}\frac{x'}{(\log x)^2}\frac{E(x)^{k_1+k_2-2}}{(k_1-1)!(k_2-1)!}.
\end{multlined}$
\end{enumerate}
\end{theorem}

\noindent\textbf{Remarques : }

\begin{enumerate}[parsep=0cm, topsep=0.2cm, label=(\roman*)]
\item Lorsque $\P$ est l'ensemble de tous les nombres premiers, on a $\N_{\P}=\N$ et donc trivialement $Q_i(n)\in\N_{\P}$ ; et $E(x)=\log_2 x+O(1)$.
\item\label{rq2th23} Pour se passer des conditions $(a_i,b_i)=1$, il suffit de diviser $a_i$ et $b_i$ par leur pgcd, en remarquant qu'alors $\omega\left(\frac{Q_i(n)}{(a_i,b_i)}\right)\in\left\{k_i,k_i-1,...,k_i-\omega((a_i,b_i))\right\}$.
\item La démonstration se généralise sans peine au cas de $\ell$ polynômes $Q_i$ de degré $1$ lorsque $\ell\geq 3$ est fixé.
\end{enumerate}

On démontre d'abord deux lemmes.

\begin{lemme}\label{lemmetenensommep}
Uniformément pour $0\leq\alpha\leq 1$, $y\geq 2$ et $\P$ un ensemble de nombres premiers, on a
\[\sum_{\substack{p\leq y \\ p\in\P}}\frac{1}{p^{1-\alpha}}=E(y)+O\left(y^{2\alpha}\right).\]
\end{lemme}

\begin{proof}
Ce résultat découle du théorème des nombres premiers par sommation d'Abel, une fois que l'on a remarqué que
\[\left\vert\sum_{\substack{p\leq y \\ p\in\P}}\frac{1}{p^{1-\alpha}}-E(y)\right\vert\leq 1+\sum_{p\leq y}\frac{1}{p^{1-\alpha}}-\frac{1}{p}.\]
\end{proof}

Étant donnés $a_1,a_2,b_1,b_2$ des entiers, on pose désormais $\Del:=a_1b_2-a_2b_1$.

\begin{lemme}\label{lemmehenriotazerty}
Soient $A, B>1$ et $0<\varepsilon<1$ fixés. Il existe deux constantes $c=c(A,\varepsilon)>0$ et $K=K(A,B,\varepsilon)>0$ vérifiant l'énoncé suivant. Soit $\P$ un ensemble de nombres premiers. Soit $x\geq 2$. Soient, pour $i\in\{1,2\}$, $Q_i(X)=a_iX+b_i$ où $1\leq a_i\leq x^A$, $0\leq b_i\leq x^A$ et $(a_i,b_i)=1$, tels que $\Del\neq 0$. Soient de plus $2\leq z_i\leq y_i\leq x$, $t_i\geq 1$ et $\xi_i:=\min\left(B, \frac{1}{4}\log y_i\right)$ pour $i\in\{1,2\}$. Soit enfin $x^{\varepsilon}<x'\leq x$. On a alors
\begin{multline*}
\left\vert\left\{n\in]x,x+x']\,:\quad\forall i\in\{1,2\},\ Q_i(n)_{y_i}>t_i,\ P^-(Q_i(n))>z_i\right\}\right\vert\\
\leq K\frac{\vert\Del\vert^K a_1a_2}{\varphi(\vert\Del\vert)^K\varphi(a_1)\varphi(a_2)}x'\prod_{i=1}^2\frac{1}{\log z_i}\exp\left(-c\xi_i\frac{\log t_i}{\log y_i}\right),
\end{multline*}
\end{lemme}

C'est une conséquence simple d'un résultat de Henriot~\cite{henriot}, qui est lui-même une version uniforme en le discriminant d'un théorème de Nair et Tenenbaum~\cite{nairtenenbaum}.

\begin{proof}
On suppose donnés les paramètres de l'énoncé. On pose
\begin{equation*}
\alpha:=\frac{\varepsilon\xi_1}{100A\log y_1},\hspace{2cm}\beta:=\frac{\varepsilon\xi_2}{100A\log y_2},
\end{equation*}
et on définit la fonction à deux variables
\[F(n,m):=\mathds{1}_{P^-\left({n}\right)>z_1}\mathds{1}_{P^-\left({m}\right)>z_2}{n}_{y_1}^{\alpha}{m}_{y_2}^{\beta},\]
qui est complètement multiplicative. Par la méthode de Rankin, le cardinal qui nous intéresse est alors
\[\leq t_1^{-\alpha}t_2^{-\beta}\sum_{x<n\leq x+x'}F(Q_1(n),Q_2(n)).\]
La fonction $F$ vérifie élémentairement
\[F(n,m)\leq\min\left(\e^{\frac{\varepsilon B}{100A}\Omega(nm)},(nm)^{\frac{\varepsilon}{400A}}\right).\]
On suppose pour simplifier qu'on est dans le cas où le polynôme $Q:=Q_1Q_2$ n'admet pas $2$ comme diviseur fixe. On peut toujours s'y ramener en séparant les cas selon que $n$ est pair ou impair. D'après~\cite[theorem~3]{henriot}, il existe alors une constante $K>0$ pouvant dépendre de $A$, $B$, et $\varepsilon$, telle que
\[\sum_{x<n\leq x+x'}F(Q_1(n),Q_2(n))\ll\frac{\vert\Del\vert^Ka_1a_2}{\varphi(\vert\Del\vert)^K\varphi(a_1)\varphi(a_2)}x'\prod_{p\leq x}\left(1-\frac{2}{p}\right)\sum_{nm\leq x}\frac{F(n,m)}{nm}.\]
Par ailleurs, d'après le Lemme~\ref{lemmetenensommep}, on a
\[\sum_{z_1<p\leq y_1}\frac{p^{\alpha}-1}{p}\ll 1,\hspace{2cm}\sum_{z_2<p\leq y_2}\frac{p^{\beta}-1}{p}\ll 1.\]
On en déduit facilement le résultat désiré. On remarque que la valeur $c=\frac{\varepsilon}{100 A}$ est admissible, indépendamment de $B$.
\end{proof}

\begin{proof}[Démonstration du Théorème~\ref{nouvelleprop}]
On ne démontre que le point~\ref{deuxiemecascri}, le premier cas en étant la variante unidimensionnelle, plus simple. On suppose donnés les paramètres de l'énoncé. On définit, pour $r\geq 1$,
\[y_r:=x^{1/2r},\]
et on pose $y_{0}=+\infty$. On trie les $n\in]x,x+x']$ selon les valeurs, pour $i\in\{1,2\}$, de $r_i\in\N$ vérifiant $Q_i(n)_{y_{r_i}}>x^{\varepsilon/3},\ Q_i(n)_{y_{r_i+1}}\leq x^{\varepsilon/3}$. On obtient alors
\[\left\{n\in]x,x+x']\,:\quad\forall i\in\{1,2\},\ \omega(Q_i(n))=k_i,\ Q_i(n)\in\N_{\P}\right\}=\bigcup_{r_1, r_2\geq 0}E_{r_1,r_2},\]
où l'on a posé, pour $r_1, r_2\geq 0$,
\begin{equation*}
E_{r_1,r_2}:=\left\{n\in]x,x+x']\,:\quad\forall i\in\{1,2\},\begin{tabular}{ll}
$\d\omega(Q_i(n))=k_i,$ & $\d Q_i(n)\in\N_{\P},$\\
$\d Q_i(n)_{y_{r_i}}>x^{\varepsilon/3},$ & $\d Q_i(n)_{y_{r_i+1}}\leq x^{\varepsilon/3}$
\end{tabular}\right\}.
\end{equation*}
On majore maintenant le cardinal des $E_{r_1,r_2}$, en commençant par  des $r_i$ non nuls, et pas trop grands. Soient $1\leq r_1,r_2\leq\sqrt{\log x}$. Étant donné un entier $N$ quelconque, on a $\omega(N_{y_{r_i+1}})\geq\omega(N)-\frac{\log N_1}{\log y_{r_1+1}}$. Ainsi, lorsque $x$ est suffisamment grand, on peut écrire
\begin{align}\label{refequne}
\vert E_{r_1,r_2}\vert&\leq\sum_{\substack{k_i-2(A+1)(r_i+1)\leq k_i'< k_i,\ \forall i\in\{1,2\} \\ k_i'\geq 0,\ \forall i\in\{1,2\}}}\sum_{\substack{d_1, d_2\leq x^{\varepsilon/3} \\ d_1d_2\in\N_{\P}\\ P^+(d_i)\leq y_{r_i+1},\ \forall i\in\{1,2\}\\ \omega(d_i)=k_i',\ \forall i\in\{1,2\}}}F_{d_1,d_2}
\end{align}
où l'on a posé
\begin{align}\label{defFd1d2}
F_{d_1,d_2}:=\left\vert\left\{n\in]x,x+x']\,:\quad\forall i\in\{1,2\},\ d_i\vert Q_i(n),\ Q_i(n)_{y_{r_i}}>x^{\varepsilon/3},\ P^-\left(\frac{Q_i(n)}{d_i}\right)>y_{r_i+1}\right\}\right\vert.
\end{align}
On note que l'inégalité stricte $k_i'<k_i$ provient du fait que si $Q_i(n)>x$ et $Q_i(n)_{y_{r_i+1}}\leq x^{\varepsilon/3}$ pour un $r_i\geq 0$, alors $Q_i(n)$ admet au moins un diviseur premier supérieur à $y_{r_i+1}$. On remarque, au vu des hypothèses, que si $(d_i,a_i)\neq 1$ pour un $i\in\{1,2\}$, alors $F_{d_1,d_2}=0$. Il en va de même si $(d_1,d_2)$ ne divise pas $\Del$. Dans le cas contraire, on paramètre $n$ sous la forme $n=\frac{d_1d_2}{(d_1,d_2)}m+\alpha$, où $\alpha:=\inf\left\{\beta\geq 0\,:\quad\forall i\in\{1,2\},\ d_i\vert a_i\beta+b_i\right\}<\frac{d_1d_2}{(d_1,d_2)}$. Pour $i\in\{1,2\}$, on note $Q_i'$ le polynôme tel que $\frac{Q_i(n)}{d_i}=Q_i'(m)$ pour ce paramétrage. On peut alors écrire
\begin{equation*}
F_{d_1,d_2}\leq F_{d_1,d_2}':=\left\vert\left\{\begin{tabular}{l}
$\d m\in\left]\frac{(x-\alpha)(d_1,d_2)}{d_1d_2},\frac{(x+x'-\alpha)(d_1,d_2)}{d_1d_2}\right]\,:$\\
$\d\forall i\in\{1,2\},\ Q_i'(m)_{y_{r_i}}>\frac{x^{\varepsilon/3}}{d_i},\ P^-(Q_i'(m))>y_{r_i+1}$
\end{tabular}\right\}\right\vert.
\end{equation*}
Soit $B>1$ une constante, que l'on choisira suffisamment grande par la suite. On pose alors, pour $i\in\{1,2\}$, $\xi_i:=\min\left(B,\frac{1}{4}\log y_{r_i}\right)$. D'après le Lemme~\ref{lemmehenriotazerty}, il existe une constante $0<c\leq 1$, indépendante de $B$, et une constante $K_0>0$ pouvant en dépendre, telles que l'on ait
\begin{equation}\label{refeqdeuxe}
F_{d_1,d_2}'\ll\frac{\vert\Del\vert^{K_0} a_1a_2}{\varphi(\vert\Del\vert)^{K_0}\varphi(a_1)\varphi(a_2)}\frac{(d_1,d_2)h(d_1)h(d_2)}{d_1d_2}\frac{x'}{(\log x)^2}\prod_{i=1}^2(r_i+1)\left(\frac{d_i}{x^{\varepsilon/3}}\right)^{\alpha_i},
\end{equation}
où l'on a posé $h(d):=\frac{d^{K_0}}{\varphi(d)^{K_0}}$, et pour $i\in\{1,2\}$,
\begin{equation}\label{defalphaiB4}
\alpha_i:=\frac{c\xi_i}{\log y_{r_i}}\leq\min\left(\frac{B}{\log y_{r_i}},\frac{1}{4}\right).
\end{equation}
La fonction $h$ est complètement sous-multiplicative et vérifie, pour $p\geq 2$ et $\nu\geq 1$, $h(p^{\nu})=1+O\left(\frac{1}{p}\right)$. En se rappelant qu'on s'est restreint au cas où $(d_1,d_2)\vert\Del$, on est alors amené à majorer
\begin{multline}\label{labelmultline}
\sum_{\substack{d_1, d_2\leq x^{\varepsilon/3},\ \forall i\in\{1,2\} \\ d_1d_2\in\N_{\P}\\ P^+(d_i)\leq y_{r_i+1},\ \forall i\in\{1,2\}\\ \omega(d_i)=k_i',\ \forall i\in\{1,2\} \\ (d_1,d_2)\vert\Del}}\frac{(d_1,d_2)h(d_1)h(d_2)}{d_1^{1-\alpha_1}d_2^{1-\alpha_2}}\\
\leq\sum_{\substack{e\vert\Del \\ P^+(e)\leq y_{r_1+1},y_{r_2+1}}}\frac{h(e)^2}{e^{1-\alpha_1-\alpha_2}}\prod_{i=1}^2\left(\sum_{\substack{d_i\in\N_{\P} \\ P^+(d_i)\leq y_{r_i+1}\\ \omega(ed_i)=k_i'}}\frac{h(d_i)}{d_i^{1-\alpha_i}}\right).
\end{multline}
On majore la dernière somme dans le cas où $\omega(e)\leq k_i'$, puisqu'elle est nulle sinon. Pour $n,m\geq 1$ des entiers, on écrit $n\vert m^{\infty}$ pour signifier que tous les facteurs premiers de $n$ divisent $m$. En décomposant $d_i$ sous la forme $d_i=f_ig_i$ où $f_i\vert e^{\infty}$ et $(g_i,e)=1$, on obtient pour $i\in\{1,2\}$,

\begin{align*}
\sum_{\substack{d_i\in\N_{\P} \\ P^+(d_i)\leq y_{r_i+1}\\ \omega(ed_i)=k_i'}}\frac{h(d_i)}{d_i^{1-\alpha_i}}&\leq\sum_{\substack{f_i\vert e^{\infty} \\ P^+(f_i)\leq y_{r_i+1}}}\frac{h(f_i)}{f_i^{1-\alpha_i}}\sum_{\substack{g_i\in\N_{\P} \\ P^+(g_i)\leq y_{r_i+1} \\ \omega(g_i)=k_i'-\omega(e)}}\frac{h(g_i)}{g_i^{1-\alpha_i}}\\
&\leq\sum_{\substack{f_i\vert e^{\infty} \\ P^+(f_i)\leq y_{r_i+1}}}\frac{h(f_i)}{f_i^{1-\alpha_i}}\frac{1}{(k_i'-\omega(e))!}\left(\sum_{\substack{p^{\nu} \\ p\in\P \\ p\leq y_{r_i+1}}}\frac{1+O\left(\frac{1}{p}\right)}{p^{\nu(1-\alpha_i)}}\right)^{k_i'-\omega(e)}\\
&\ll\sum_{\substack{f_i\vert e^{\infty} \\ P^+(f_i)\leq y_{r_i+1}}}\frac{h(f_i)}{f_i^{1-\alpha_i}}\frac{E(x)^{k_i'-\omega(e)}}{(k_i'-\omega(e))!}\\
&\ll\frac{E(x)^{k_i'}}{k_i'!}A^{\omega(e)}\sum_{\substack{f_i\vert e^{\infty} \\ P^+(f_i)\leq y_{r_i+1}}}\frac{h(f_i)}{f_i^{1-\alpha_i}}\\
&\ll\frac{E(x)^{k_i'}}{k_i'!}A^{\omega(e)}\left(\frac{e}{\varphi(e)}\right)^{\e^B}.
\end{align*}
La troisième inégalité est obtenue en utilisant le Lemme~\ref{lemmetenensommep} et~$(\ref{defalphaiB4})$, puis la croissance de $x\mapsto E(x)$ et le fait que $k_i'\leq AE(x)$. La dernière inégalité est obtenue en passant au produit eulérien et en utilisant de nouveau~$(\ref{defalphaiB4})$. En reportant dans~$(\ref{labelmultline})$ et en utilisant encore~$(\ref{defalphaiB4})$, on obtient
\[\sum_{\substack{d_1, d_2\leq x^{\varepsilon/3},\ \forall i\in\{1,2\} \\ d_1d_2\in\N_{\P}\\ P^+(d_i)\leq y_{r_i+1},\ \forall i\in\{1,2\}\\ \omega(d_i)=k_i',\ \forall i\in\{1,2\} \\ (d_1,d_2)\vert\Del}}\frac{(d_1,d_2)h(d_1)h(d_2)}{d_1^{1-\alpha_1}d_2^{1-\alpha_2}}\ll\frac{E(x)^{k_1'+k_2'}}{k_1'!k_2'!}\left(\frac{\vert\Del\vert}{\varphi(\vert\Del\vert )}\right)^{A\e^{2B}}\]
Finalement, en reportant~$(\ref{refeqdeuxe})$ dans~$(\ref{refequne})$, puisque $k_i'\leq k_i-1\leq AE(x)$, on obtient
\begin{multline}\label{qsjdhskqjdh}
\vert E_{r_1,r_2}\vert\ll\frac{\vert\Del\vert^{K_0+A\e^{2B}}a_1a_2}{\varphi(\vert\Del\vert)^{K_0+A\e^{2B}}\varphi(a_1)\varphi(a_2)}\frac{x'}{(\log x)^2}\frac{E(x)^{k_1+k_2-2}}{(k_1-1)!(k_2-1)!}\\
\times\prod_{i=1}^2A^{2(A+1)(r_i+1)}(r_i+1)^2\exp\left(-\frac{c\varepsilon}{12}\min\left(Br_i,\log x\right)\right).
\end{multline}
Lorsque $B$ est suffisamment grand -- on rappelle que la constante $c$ est indépendante de $B$ --, on en déduit
\[\sum_{1\leq r_1,r_2\leq\sqrt{\log x}}\vert E_{r_1,r_2}\vert\ll\frac{\vert\Del\vert^{K_0+A\e^{2B}}a_1a_2}{\varphi(\vert\Del\vert)^{K_0+A\e^{2B}}\varphi(a_1)\varphi(a_2)}\frac{x'}{(\log x)^2}\frac{E(x)^{k_1+k_2-2}}{(k_1-1)!(k_2-1)!}.\]
Le cas des $E_{r_1,r_2}$ où l'un des $r_i$ vaut $0$ se traite de manière analogue. Il suffit, dans la définition~$(\ref{defFd1d2})$ de $F_{d_1,d_2}$, d'ignorer la condition $Q_i(n)_{y_0}>x^{\varepsilon/3}$, ce qui revient, dans~$(\ref{qsjdhskqjdh})$, à remplacer l'exponentielle par $1$. En posant $r:=\sqrt{\log x}$, on remarque enfin que l'on a
\[\bigcup_{r_1\text{ ou }r_2\geq r}E_{r_1,r_2}\subset E^{(1)}\cup E^{(2)},\]
où l'on a posé, pour $i\in\{1,2\}$,
\[E^{(i)}:=\left\{n\in]x,x+x']\,:\quad Q_i(n)_{y_r}>x^{\varepsilon/3}\right\}.\]
En utilisant la méthode de Rankin comme précédemment et~\cite[theorem~3]{henriot}, on obtient alors la majoration adéquate pour les $\vert E^{(i)}\vert$, puisque $\frac{1}{(\log x)^2}\frac{E(x)^{k_1+k_2-2}}{(k_1-1)!(k_2-1)!}\gg\exp\left(-(\log x)^{1/4}\right)$, disons.

\end{proof}

\subsection{Application du Théorème~\ref{th1}}

On vérifie maintenant, à l'aide des résultats précédents, que l'on peut appliquer le Théorème~\ref{th1} pour démontrer le Théorème~\ref{corollaire}.

\begin{proof}[Démonstration du Théorème~$\ref{corollaire}$]
Le résultat est directement vrai lorsque $X<h\leq \delta_k(X)^{-1}X$ d'après le Lemme~$\ref{hyp1application}$. Soit $R>1$ fixé. Soient $X\geq 3$ et $R^{-1}\log_2 X\leq k\leq R\log_2 X$. Soient $\psi:\R_+\rightarrow\R_+$ tendant vers l'infini en l'infini et $1\leq\psi(X)\leq\delta_k(X)h\leq X$. On définit alors les paramètres suivants pour appliquer le Théorème~$\ref{th1}$. On fixe $r=\frac{1}{2R}$ et $0<\varepsilon\leq\frac{1}{100}$ et on prend $K=(\log_2 X)^2$, de sorte que $J\ll\log_3 X$, et $P_{\infty},Q_{\infty}$ vérifiant~$(\ref{defS2})$. On pose $\mathcal{A}=\mathcal{E}_k$ et on prend $T_0:=\left(\log X\right)^{100+R^2}$ et $f:=\left(\frac{\mathrm{Id}}{\varphi}\right)^{K_1}$ où $K_1$ est une constante à choisir suffisamment grande. On prend $\theta=1$ et $\delta=\delta_k$. Enfin, on définit
\begin{equation}\label{defX0}
P_1:=\exp\left(\frac{\log Q_1}{1+\min\left(\sqrt{\log\psi(X)},\log_3 Q_1-1\right)}\right),
\end{equation}
qui satisfait à~$(\ref{defS2})$. On vérifie désormais les hypothèses du Théorème~$\ref{th1}$. La première est vérifiée d'après le Lemme~$\ref{hyp1application}$. Soit désormais $n\in\left(\mathcal{A}\smallsetminus\mathcal{S}\right)\cap\left[X,2X\right]$. Au moins une des trois possibilités suivantes est vérifiée pour au moins un indice $j\in[1,J]\cup\{\infty\}$ :
\begin{itemize}
\item aucun $p\in\left]P_j,Q_j\right]$ ne divise $n$,
\item il existe un $p\in\left]P_j,Q_j\right]$ tel que $p^2\vert n$,
\item il existe $p,q\in\left]P_j,Q_j\right]$ tels que $pq\vert n$ et $\left\vert p-q\right\vert\ll\frac{p}{H_j}$.
\end{itemize}
D'après le Lemme~$\ref{hyp2application}$, le nombre d'entiers vérifiant la première possibilité pour au moins un indice $j$ est
\[\ll\delta_k(X)X\sum_{j\in\N\cup\{\infty\}}\left(\frac{\log P_j}{\log Q_j}\right)^r.\]
Pour la deuxième possibilité, il y a moins de
\[\ll\sum_{j\in\N\cup\{\infty\}}\sum_{P_j<p\leq Q_j}\left\vert\left\{n\sim\frac{X}{p^2}\,:\quad\omega(p^2n)=k\right\}\right\vert\ll\delta_k(X)X\sum_{j\in\N\cup\{\infty\}}\frac{1}{P_j}\]
entiers, la deuxième inégalité étant obtenue en observant d'une part que $Q_{\infty}=X^{o(1)}$, et d'autre part que $\delta_{k-1}(X)\ll\delta_k(X)$, $R$ étant fixé. Enfin, avec les même remarques, le nombre d'entiers de $\left(\mathcal{A}\smallsetminus\mathcal{S}\right)\cap\left[X,2X\right]$ vérifiant le dernier point pour au moins un $j$ est
\[\ll\sum_{j\in\N\cup\{\infty\}}\sum_{P_j<p\leq Q_j}\sum_{\substack{P_j<q\leq Q_j \\ \left\vert p-q\right\vert\ll\frac{p}{H_j}}}\left\vert\left\{n\sim\frac{X}{pq}\,:\quad\omega(pqn)=k\right\}\right\vert\ll\delta_k(X)X\sum_{j\in\N\cup\{\infty\}}\frac{\log\left(1+\frac{\log Q_j}{\log P_j}\right)}{H_j}.\]
L'hypothèse~\ref{hyp2} est donc vérifiée. La troisième hypothèse est directement vérifiée d'après le Théorème~\ref{nouvelleprop} lorsque $K_1$ est suffisamment grand. Ensuite, pour tous nombres premiers $p,q$ et tout entier $m$ tels que $q^2\nmid mq$ et $p\nmid m$, on a $\omega(mp)=\omega(mq)$. L'hypothèse~\ref{hyp4} est donc vérifiée. En remarquant que $\mathcal{A}'\subset\mathcal{E}_{k-1}$, l'hypothèse~\ref{hyp5} est aussi facilement vérifiée. Il en va de même pour l'hypothèse~\ref{hyp6} puisque $\delta_{k+\ell}(X)\ll r^{-\ell}\delta_k(X)$. Enfin, les hypothèses~\ref{hyp7} et~\ref{hyp8} sont aisément vérifiées. On peut donc appliquer le Théorème~$\ref{th1}$ et en tirer la conclusion désirée.
\end{proof}

\section{Entiers friables}\label{applifriall}

On démontre ici les Théorèmes~$\ref{thfriables}$ et~$\ref{tsintervallesfriables}$. Il semble très difficile, à l'heure actuelle, de majorer efficacement $\left\vert\left\{ n\sim x\,:\quad P^+(n(n+1))\leq y\right\}\right\vert$ pour un paramètre $y$ donné, même relativement proche de $x$. C'est une différence fondamentale avec le cas traité précédemment des entiers ayant un nombre donné $k$ de facteurs premiers. On peut quand même tirer une information du Théorème~$\ref{th1}$, en prenant la perte $\theta$ maximale, afin de trivialiser quelques hypothèses. Dans l'optique de vérifier les autres, on énonce trois lemmes.

\begin{lemme}[Hildebrand,~\cite{hildebrandfriables}]\label{hild}
Soit $\psi:\R_+\rightarrow\R_+$ une fonction tendant vers l'infini en l'infini. Soit $\varepsilon>0$ fixé. Pour $x\geq 3, 1\leq u\leq\frac{\log x}{(\log_2 x)^{3/5+\varepsilon}}$ et $1\leq z\leq x^{5/12u}$, on a
\[\Psi\left(x+\frac{x}{z},x^{1/u}\right)-\Psi\left(x,x^{1/u}\right)=\rho(u)\frac{x}{z}\left(1+o(1)\right)\]
lorsque $x$ tend vers l'infini.
\end{lemme}

On a l'estimation suivante de la fonction de Dickman~\cite{HildTenensurvolfriables}.

\begin{lemme}\label{explicitrho}
La fonction $\rho:\R\rightarrow\R$ de Dickman, définie comme l'unique solution continue sur $\R_+$ du système
\begin{equation*}
\left\{\begin{tabular}{ll}
$u\rho'(u)=-\rho(u-1)$ & pour $u>1$,\\
$\rho(u)=1$ & pour $0\leq u\leq 1$,\\
$\rho(u)=0$ & pour $u<0$,
\end{tabular}\right.
\end{equation*}
vérifie, pour $u\geq 1$,
\begin{equation*}\label{eqrho}
\rho(u)=\exp\left(-u\left(\log u+\log_2\left(u+2\right)-1+O\left(\frac{\log_2\left(u+2\right)}{\log\left(u+2\right)}\right)\right)\right).
\end{equation*}
\end{lemme}

\begin{lemme}\label{sommefriables}
Soit $0<\varepsilon<\frac{1}{6}$ fixé. Uniformément pour $x\geq 3$, $1\leq u\leq (\log x)^{1/6-\varepsilon}$, $1\leq z\leq \exp\left((\log x)^{1/6-\varepsilon/2}\right)$ et $2\leq P\leq Q\leq \exp\left((\log x)^{5/6+\varepsilon/3}\right)$, on a
\begin{multline*}
\left\vert\left\{n\in\left]x, x+\frac{x}{z}\right]\,:\quad P^+(n)\leq x^{1/u},\ \left(n,\prod_{P<p\leq Q}p\right)=1\right\}\right\vert\\
=\left(1+O\left(\frac{1}{u}\right)\right)\prod_{P<p\leq Q}\left(1-\frac{1}{p}\right)\rho(u)\frac{x}{z}.
\end{multline*}
\end{lemme}

Ce lemme est une extension aux petits intervalles d'un résultat beaucoup plus général que La Bretèche \& Tenenbaum~\cite{BretecheTenenbaumfriables} obtiennent dans leur article lorsque $z=1$. 

\begin{proof}
On suppose donnés les paramètres de l'énoncé. Comme de coutume, pour $k\geq 1$, on définit $\rho^{(k)}$ sur $\R$ tout entier par continuité à droite. La démonstration est une conséquence rapide des méthodes développées dans~\cite{BretecheTenenbaumfriables}. On en reprend les notations. Notamment, pour $x, y>0$ des réels et $m$ un entier, on définit
\begin{align*}
\Psi_m(x,y)&:=\sum_{\substack{n\leq x \\ P^+(n)\leq y\\ (n,m)=1}}1,\\
\Lambda_m(x,y)&:=\left\{\begin{tabular}{ll}
$\d x\int_{0-}^{\infty}\rho(u-v)\mathrm{d}R_m(y^v)$ & $\d(x\in\R\smallsetminus\Z),$ \\ $\d\Lambda_m(x+,y)$ & $\d(x\in\Z),$
\end{tabular}\right.\\
R_m(x)&:=\frac{1}{x}\sum_{\substack{n\leq x \\ (n,m)=1}}1-\frac{\varphi(m)}{m}.
\end{align*}
On pose ici $m=\prod_{P<p\leq Q}p$.  On commence par utiliser~$(4\cdot1)$ de~\cite{BretecheTenenbaumfriables}, afin d'obtenir
\begin{multline*}
\Psi_m\left(x+\frac{x}{z},x^{1/u}\right)-\Psi_m\left(x,x^{1/u}\right)=\Lambda_m\left(x+\frac{x}{z},x^{1/u}\right)-\Lambda_m\left(x,x^{1/u}\right)\\
+O\left(\rho(u)x\exp\left(-\sqrt{\log x}\right)\right).
\end{multline*}
À partir de~$(4\cdot13)$ de~\cite{BretecheTenenbaumfriables}, un calcul élémentaire montre qu'il suffit alors de prouver
\begin{equation}\label{eqvmdiff}
V_m\left(x+\frac{x}{z},x^{1/u}\right)-V_m\left(x,x^{1/u}\right)\ll\frac{\rho(u)\varphi(m)}{uzm},
\end{equation}
où l'on a posé
\begin{equation*}
V_m(x,y):=R_m(x)+\int_0^{\infty}\left\{\rho'(u-v)-\rho'(u)\e^{r(u)v}\right\}R_m(y^v)\mathrm{d}v,
\end{equation*}
et, pour $u>0$,
\[r(u):=\frac{-\rho'(u)}{\rho(u)}=\frac{\rho(u-1)}{u\rho(u)}.\]
Puisque $P^+(m)\leq\exp\left((\log x)^{5/6+\varepsilon/3}\right)$, l'inégalité~$(3\cdot41)$ de~\cite[lemme~$3.10$]{BretecheTenenbaumfriables} fournit
\[\left\vert R_m\left(x+\frac{x}{z}\right)\right\vert+\left\vert R_m(x)\right\vert\ll\frac{\varphi(m)}{m}\exp\left(-(\log x)^{1/6-\varepsilon/3}\right),\]
qui est convenable. On définit $u'$ tel que
\[\left(x+\frac{x}{z}\right)^{1/u'}=x^{1/u},\]
de sorte que $u\leq u'\leq u+O\left(\frac{u}{z\log x}\right)$, et on pose $w(u,v):=\rho'(u-v)-\rho'(u)\e^{r(u)v}$. On écrit alors le membre de gauche de~$(\ref{eqvmdiff})$ sous la forme $I_1+I_2+I_3+R$ où l'on a posé
\begin{align*}
I_1&:=\int_0^{u-2}\left\{w(u',v)-w(u,v)\right\}R_m(x^{v/u})\mathds{1}_{v\geq 0}\mathrm{d}v,\\
I_2&:=\int_{u-2}^{u'-1}\left\{\rho'(u'-v)-\rho'(u-v)\right\}R_m(x^{v/u})\mathds{1}_{v\geq 0}\mathrm{d}v,\\
I_3&:=\int_{u-2}^{\infty}\left\{\rho'(u')\e^{r(u')v}-\rho'(u)\e^{r(u)v}\right\}R_m(x^{v/u})\mathds{1}_{v\geq 0}\mathrm{d}v,\\
R&:=R_m\left(x+\frac{x}{z}\right)- R_m(x).
\end{align*}
La majoration de $I_2$ ne pose aucun problème à partir des lemmes~$3.9$ et~$3.10$ de~\cite{BretecheTenenbaumfriables}. On majore désormais $I_1$. Montrons, pour $u\geq 1$ et $0\leq v\leq u-2$, la majoration
\begin{align}\label{majderivee}
\frac{\partial w}{\partial u}(u,v)&=\rho''(u-v)-\rho''(u)\e^{r(u)v}-vr'(u)\rho'(u)\e^{r(u)v}\nonumber\\
&\ll\log(u+1)\frac{\rho(u)v\e^{r(u)v}}{u}\left(1+v\log(u+1)\right).
\end{align}
Soit
\[\lambda(u):=\rho''(u)+\frac{\rho'(u)}{u}=-\frac{\rho'(u-1)}{u}=\frac{r(u-1)\rho(u-1)}{u}.\]
Le lemme~6.1 de~\cite{repstatFouTen} ainsi que sa démonstration fournissent
\[\lambda(u-v)=\lambda(u)\e^{r(u)v}\left\{1+O\left(\frac{v}{u\log(u+v)}+\frac{v^2}{u}\right)\right\}.\]
Le même lemme, combiné avec~\cite[corollaire~\chrom{3}.5.14]{tenenbaum}, nous permet alors d'écrire
\[\rho''(u-v)-\rho''(u)\e^{r(u)v}\ll\log(u+1)\frac{\rho(u)v\e^{r(u)v}}{u}\left(1+v\log(u+1)\right).\]f
L'inégalité~$(\ref{majderivee})$ est alors facilement obtenue une fois connue la majoration $r'(u)\ll\frac{1}{u}$, qui correspond à~$(6.9)$ de~\cite[lemme~6.1]{repstatFouTen}. De manière analogue à~\cite[$(4\cdot19)$]{BretecheTenenbaumfriables}, la majoration adéquate pour $I_1$ découle alors directement du théorème des accroissements finis, puisque $\log u\leq\log (x^{1/u})$.

La majoration de $I_3$ s'effectue de façon analogue, via les accroissements finis, en utilisant le travail de majoration qui a été effectué dans~\cite{BretecheTenenbaumfriables} pour obtenir~$(4\cdot18)$. Ces méthodes n'étant pas au centre du présent article, nous laissons les détails au lecteur intéressé.
\end{proof}

On montre maintenant le Théorème~$\ref{thfriables}$.

\begin{proof}[Démonstration du Théorème~$\ref{thfriables}$]
Soit $0<\varepsilon<\frac{1}{100}$ fixé. Soient $X\geq 3$ et $1\leq u\leq (\log X)^{1/6-\varepsilon}$. On applique le Théorème~$\ref{th1}$ avec les paramètres suivants. On pose $r=1$, $K=\sqrt{\log X}$, $P_{\infty}=\exp\left((\log X)^{5/6-\varepsilon/6}\right)$ et $Q_{\infty}=\exp\left((\log X)^{5/6}\right)$. On prend $\mathcal{A}$ l'ensemble des entiers $X^{1/u}$-friables et $\delta(X):=\rho(u)$. On remarque que l'on a $\delta(X)\gg\exp\left(-(\log X)^{1/6-5\varepsilon/6}\right)\gg P_{\infty}^{-(\log X)^{-2/3-\varepsilon/2}}$ et on pose $T_0:=\exp\left((\log X)^{1/6}\right)$. On pose $\theta(X):=\delta(X)^{-1}$ et $f:=1$. Soient $\psi:\R_+\rightarrow\R_+$ une fonction tendant vers l'infini en l'infini et $\left(1+\rho(u)^{-1}\right)^{\psi(X)}\leq h\leq X$. Enfin, on pose
\[P_1:=\exp\left(\frac{\log Q_1}{1+\min\left(\sqrt{\psi(X)},\log_3 Q_1-1\right)}\right),\]
qui satisfait à~$(\ref{defS2})$. On vérifie désormais les hypothèses du Théorème~$\ref{th1}$. Le premier point découle aisément du Lemme~$\ref{hild}$. Pour vérifier le deuxième point, on procède exactement comme pour le Théorème~$\ref{corollaire}$ en utilisant les Lemmes~$\ref{hild}$ et~$\ref{sommefriables}$. Puisque $\theta=\delta^{-1}$, le troisième point est trivialement vérifié. L'hypothèse~$\ref{hyp4}$ est vérifiée puisque $\mathcal{A}$ est représenté par une fonction complètement multiplicative, valant $1$ sur tous les premiers $p\leq Q_{\infty}$. Les points~\ref{hyp5} et~\ref{hyp6} ne posent aucun problème puisque $\mathcal{A}'\subset\mathcal{A}$. Enfin, les hypothèses~\ref{hyp7} et~\ref{hyp8} se vérifient de manière élémentaire.
\end{proof}

\begin{proof}[Démonstration du Théorème~$\ref{tsintervallesfriables}$]
Soit $0<\varepsilon<\frac{1}{1000}$ fixé. Soient $x\geq 1$ et $1\leq u\leq (\log x)^{1/6-\varepsilon}$. Soient $\rho\left(\frac{u}{2}\right)^{-6-100\varepsilon}\leq h_1\leq\sqrt{x}$, puis $T_0:=\rho\left(\frac{u}{2}\right)^{-20}(\log x)^{10}$ et $h_2:=\frac{\sqrt{x}}{T_0^3}$. Comme dans~\cite{matoradzi}, on étudie, pour $i\in\{1,2\}$,
\[\Sigma(h_i):=\frac{1}{h_i\sqrt{x}}\sum_{\substack{x-h_i\sqrt{x}<n_1n_2\leq x+h_i\sqrt{x} \\ n_1\sim\sqrt{x} \\ n_1,n_2\in\mathcal{S} \\ P^+(n_1n_2)\leq x^{1/u}}}1\]
où $\mathcal{S}$ est un ensemble que l'on définit plus tard. On traite d'abord le cas $h_1\leq h_2$. En procédant de la même manière que dans~\cite[proof of theorem~4]{matoradzi} --- en prenant leur paramètre $\delta$ égal à $1$ ---, on obtient
\begin{equation*}\label{eqSigma}
2\Sigma(2h_1)\gg\Sigma(h_2)+O\left(\frac{1}{T_0}+I\right),
\end{equation*}
où
\[I:=\int_{T_0}^{\frac{\sqrt{x}}{h_1}}\left\vert M_1(1+it)M_2(1+it)\right\vert\mathrm{d}t+\frac{\sqrt{x}}{h_1}\max_{T>\frac{\sqrt{x}}{h_1}}\frac{1}{T}\int_T^{2T}\left\vert M_1(1+it)M_2(1+it)\right\vert\mathrm{d}t\]
avec les notations
\[M_1(s):=\sum_{\substack{n\sim\sqrt{x} \\ n\in\mathcal{S} \\ P^+(n)\leq x^{1/u}}}\frac{1}{n^s},\hspace{1cm}M_2(s):=\sum_{\substack{\frac{\sqrt{x}}{2}<n\leq 2\sqrt{x} \\ n\in\mathcal{S} \\ P^+(n)\leq x^{1/u}}}\frac{1}{n^s}\hspace{1cm}(s\in\C).\]
En utilisant l'inégalité de Cauchy-Schwarz, on est ramené à majorer, pour $i\in\{1,2\}$ et $T\geq T_0$,
\[\int_{T_0}^{T}\left\vert M_i(1+it)\right\vert^2\mathrm{d}t.\]
On pose $X:=\sqrt{x}$. On définit maintenant $\mathcal{S}$ comme dans la section~$\ref{section3}$ pour les paramètres $r=1$, $\varepsilon$ donné ci-dessus, $\delta=\rho\left(\frac{u}{2}\right),\theta=\delta^{-1}, h=h_1, K=(\log X)^{1/6-2\varepsilon/3}, P_1=Q_1^{1-\varepsilon},  P_{\infty}=\exp\left((\log X)^{5/6-\varepsilon/6}\right)$ et $Q_{\infty}=\exp\left((\log X)^{5/6}\right)$ ; à une différence près, en posant, pour $1\leq j\leq J$, $H_j:=j^2P_1^{1/6-\varepsilon}(\log Q_1)^{-1/3}$. On remarque alors, comme dans la démonstration du théorème précédent, que les hypothèses~\ref{hyp1} à~\ref{hyp6} sont toutes vérifiées pour $u$ assez grand, l'hypothèse~\ref{hyp3} étant même vérifiée pour $\frac{X}{H_J}\ll y\leq X$. Cela garantit qu'une démonstration identique à celle de la Proposition~$\ref{majintegrale}$ fournisse, pour $i\in\{1,2\}$,
\[\int_{T_0}^{T}\left\vert M_i(1+it)\right\vert^2\mathrm{d}t\ll\left(\frac{TQ_1}{\delta(X)X}+\theta(X)\right)\delta(X)^2P_1^{-1/6+2\varepsilon},\]
pour $T_0<T\leq X$. Ainsi, en utilisant directement le Lemme~$\ref{MVT}$ lorsque $T>X$, on obtient
\[I\ll\rho\left(\frac{u}{2}\right)^{2+\varepsilon}\]
lorsque $u$ est assez grand. Par ailleurs, avec les Lemmes~\ref{hild},~\ref{explicitrho} et~\ref{sommefriables}, on trouve
\[\Sigma(h_2)\gg\rho\left(\frac{u}{2}\right)^2.\]
Finalement, on obtient le résultat désiré lorsque $u$ est assez grand en observant que l'on a
\begin{align*}
\rho\left(\frac{u}{2}\right)^{2}h_1\sqrt{x}&\ll\sum_{\substack{x-2h_1\sqrt{x}<n_1n_2\leq x+2h_1\sqrt{x} \\ n_1\sim\sqrt{x} \\ P^+(n_1n_2)\leq x^{1/u}}}1\\
&\ll\left(\sum_{\substack{\left\vert n-x\right\vert\leq 2h_1\sqrt{x} \\ P^+(n)\leq x^{1/u}}}1\right)^{1/2}\left(\sum_{\left\vert n-x\right\vert\leq 2h_1\sqrt{x}}\left(\sum_{n_1n_2=n} 1\right)^2\right)^{1/2}\\
&\ll\left(\sum_{\substack{\left\vert n-x\right\vert\leq 2h_1\sqrt{x} \\ P^+(n)\leq x^{1/u}}}1\right)^{1/2}\left(h_1\sqrt{x}(\log x)^3\right)^{1/2},
\end{align*}
et $\rho\left(\frac{u}{2}\right)^2\gg\rho(u)2^{u+o(u)}$ lorsque $u$ tend vers l'infini. La dernière inégalité ci-dessus est obtenue par application directe de~\cite[theorem~3]{henriot}. En effet, ce théorème fournit, pour $x^{\varepsilon}<x'\leq x$, l'inégalité
\[\sum_{x<n\leq x+x'}\tau(n)^2\ll x'\prod_{p\leq x}\left(1-\frac{1}{p}\right)\sum_{n\leq x}\frac{\tau(n)^2}{n}.\]
Le cas où $h_2\leq h_1\leq\sqrt{x}$ se traite directement en observant que pour $\mathcal{S}=\N$, le Lemme~\ref{hild} fournit directement $\Sigma(h_1)\gg\rho\left(\frac{u}{2}\right)^{2}$.

\end{proof}

\section{Entiers avec un petit nombre de facteurs premiers}\label{kborne}

De manière informelle, un des ingrédients clefs du Théorème~$\ref{th1}$ est le suivant. Presque tous les entiers que l'on étudie ont au moins $1$ facteur premier dans chacun des $J$ intervalles $]P_j,Q_j]$. Cela met à notre disposition $J$ factorisations possibles pour le polynôme de Dirichlet $B$, ce qui nous permet de majorer 
\[\int_{T_0}^T\vert B(1+it)\vert^2\mathrm{d}t.\]
Pour $t\in\left[T_0,T\right]$, on se sert de l'indice $j$ fournissant la factorisation la plus avantageuse. C'est la factorisation correspondant à l'indice $j=1$ qui est la plus difficile à majorer, et qui impose à $h$ d'être suffisamment grand, essentiellement de taille $\delta(X)^{-1}Q_1$. Plus $Q_1$ est petit, plus $h$ peut être petit. Par ailleurs, pour que la factorisation correspondant à un certain $j\geq 2$ soit efficace, il faut que $P_{j-1}$ et $Q_{j-1}$ ne soient pas \og trop éloignés\fg\ de $P_j$ et $Q_j$. Ainsi, plus on veut réduire la taille de $Q_1$, et donc de $h$, plus $J$ doit être grand. Dans le cadre du Théorème~$\ref{th1}$, ce n'est pas contraignant, grâce à l'hypothèse~\ref{hyp2}. On a déjà vu que lorsque l'on s'intéressait aux entiers ayant $k\asymp\log_2 X$ facteurs premiers, cette hypothèse était vérifiée (\emph{cf.} Lemme~\ref{hyp2application}). Cela revient à dire, lorsque $\log P=o\left(\log Q\right)$, que presque tous les entiers ayant $k$ facteurs premiers en ont au moins $1$ dans l'intervalle $]P,Q]$. Ainsi, il est suffisant de s'intéresser à ces entiers seulement, ce qui revient à travailler dans $\mathcal{A}\cap\mathcal{S}$. 

Lorsque $k$ est trop petit, par exemple borné, cela n'est évidemment plus vrai. D'une part, on ne peut plus espérer d'un entier avec $k$ facteurs premiers d'en avoir $J$ localisés dans certains intervalles si $J>k$. D'autre part, il n'y a approximativement qu'une proportion $\frac{\log_2 Q-\log_2 P}{\log_2 X}$ des entiers avec $2\leq k\ll 1$ facteurs premiers qui en ont un dans l'intervalle $]P,Q]$. La restriction à $\mathcal{S}$ n'est donc plus innocente dans ce cas-là. Il paraît alors difficile, par cette méthode, de comparer $\frac{1}{h}\left(\pi_k(x+h)-\pi_k(x)\right)$ à $\frac{1}{X}\left(\pi_k(2X)-\pi_k(X)\right)$, pour $k$ petit. On peut cependant minorer la première quantité en comparant
\[\frac{1}{h}\sum_{\substack{x<n\leq x+h \\ n\in\mathcal{E}_k\cap\mathcal{S}}}1\hspace{1cm}\text{\ à\ }\hspace{1cm}\frac{1}{X}\sum_{\substack{n\sim X \\ n\in\mathcal{E}_k\cap\mathcal{S}}}1.\]
Par la remarque précédente, $\mathcal{E}_k\cap\mathcal{S}$ est de densité nulle dans $\mathcal{E}_k$. De manière approximative, dès que $Q$ n'est pas proche de $X$, le fait d'imposer un facteur premier dans $]P,Q]$ fait perdre un facteur $\frac{\log_2 X}{k}$ sur la densité. Il faut donc choisir $J$ minimal de sorte que $Q_1$ soit plus petit que disons $\log_2 X$. On rappelle la décomposition $[T_0,T]=\bigcup_{1\leq j\leq J}\mathcal{T}_j\cup\mathcal{U}$. La majoration de l'intégrale sur $\mathcal{U}$ est aisée dès que l'on a, approximativement, $\log_2 X=o\left(\log P_J\right)$, et la majoration sur $\mathcal{T}_j$ pour $j\geq 2$ dès que $\log_2 Q_j=o\left(\log P_{j-1}\right)$. Cela implique approximativement $Q_1\geq(\log_J X)^{C}$. Par ce raisonnement, on espère minorer $\frac{1}{h}\left(\pi_k(x+h)-\pi_k(x)\right)$ dès que $\delta_k(X)h\geq\left(\frac{\log_2 X}{k}\right)^J(\log_J X)^C$. Notre intérêt portant ici sur les petites valeurs de $k$, on choisit le paramètre $J$ égal à $3$. Comme mentionné en introduction, pour des valeurs de $k$ dépassant $\frac{\log_2 X}{\log_3 X}$, il est possible d'améliorer légèrement le Théorème~\ref{petitsk}, en prenant le paramètre $J$ plus grand.

Comme on vient de le voir, la démonstration ressemble très fortement à celle du Théorème~$\ref{th1}$, mais avec un choix légèrement différent pour~$\mathcal{S}$. On se permet alors d'omettre certains calculs, qui ressemblent fortement à ceux déjà effectués dans la section~$\ref{sectionth1}$.

\subsection{Redéfinition de $\mathcal{S}$}\label{sssec}

Soit $\varepsilon>0$ fixé. Pour $X$ suffisamment grand, on définit
\begin{equation*}\label{spetitsk}
\left\{\begin{tabular}{llr}
$P_1:=(\log_3 X)^{2+100\varepsilon}$, & $Q_1:=P_1^{1+\varepsilon}$,\\
$P_2:=(\log_2 X)^{\varepsilon^{-1}}$, & $Q_2:=P_2^{1+\varepsilon}$,\\
$P_3:=(\log X)^{\varepsilon^{-1}}$, & $Q_3:=\exp\left((\log_2 X)^2\right)$,\\
$P_{\infty}:=\exp\left((\log X)^{2/3+\varepsilon}\right)$, & $Q_{\infty}:=\exp\left((\log X)^{2/3+2\varepsilon}\right)$,\\
$H_i:=(\log_2 Q_i)^2$ & & $(j\in[1,3]\cup\{\infty\})$.
\end{tabular}\right.
\end{equation*}
On définit alors $\mathcal{S}$ comme dans la section~$\ref{section3}$. 
\begin{equation*}
\left\{\begin{tabular}{lr}
$\mathcal{I}_j:=\left[\left\lfloor H_j\log P_j\right\rfloor,H_j\log Q_j\right]$, & $(j\in[1,3]\cup\{\infty\})$\\
$\d\mathcal{S}_j:=\bigcap_{v\in\mathcal{I}_j}\bigcap_{\substack{\e^{v/H_j}<p,q\leq\e^{(v+1)/H_j} \\ P_j<p,q\leq Q_j}}\left\{n\geq 1\,:\quad\mu_{]P_j,Q_j]}^2\omega_{]P_j,Q_j]}(n)\geq 1,\ pq\nmid n\right\}$, & $(j\in[1,3]\cup\{\infty\})$
\end{tabular}\right.
\end{equation*}
et
\begin{equation*}
\mathcal{S}:=\bigcap_{j=1}^3\mathcal{S}_j\cap\mathcal{S}_{\infty}.
\end{equation*}
Pour la suite on pose
\begin{equation}\label{alphapetitsk}
\left\{\begin{tabular}{l}
$\alpha_1:=\varepsilon$,\\
$\alpha_2:=\frac{1}{4}-4\varepsilon$,\\
$\alpha_3:=\frac{1}{4}-2\varepsilon$.
\end{tabular}\right.
\end{equation}
Pour $5\leq k\leq\log_2 X$, on rappelle la définition~$(\ref{deffkdjfhkjdfghdk})$ de~$F_k$. On note que si l'on avait choisi $Q_3$ égal à $P_3^{1+\varepsilon}$, l'énoncé du Théorème~\ref{petitsk} serait toujours valable mais en remplaçant $F_k(X)$ par $\frac{(\log_2 X)^3}{k^3}$. Le choix que l'on fait effectivement permet de gagner un facteur $\log_3 X$ lorsque $k=o\left(\frac{\log_2 X}{\log_3 X}\right)$.

\subsection{Démonstration du Théorème~$\ref{petitsk}$}

On commence par estimer la densité de $\mathcal{E}_k\cap\mathcal{S}$.

\begin{lemme}\label{sommepetitsk}
On conserve les notations de la sous-section~$\ref{sssec}$, et on pose $T_0:=(\log X)^{10}$ et $y_0:=\frac{X}{T_0^3}$. On a alors, pour tout $x\sim X$,
\[\sum_{\substack{x<n\leq x+y_0 \\ n\in\mathcal{E}_k\cap\mathcal{S}}}1\asymp\frac{\delta_k(X)}{F_k(X)}y_0.\]
\end{lemme}

\begin{proof}
Tout entier $n\in\mathcal{E}_k\cap\mathcal{S}$ se décompose de manière unique sous la forme $m_1m_2m_3r$ où, pour $i\in[1,3]$, l'entier $m_i$ correspond au produit des facteurs premiers de $n$ appartenant à $]P_i,Q_i]$. En notant $m:=m_1m_2m_3$, on a alors
\[\sum_{\substack{x<n\leq x+y_0 \\ n\in\mathcal{E}_k\cap\mathcal{S}}}1=\sum_{\substack{m_i, \forall i\in[1,3] \\ p\vert m_i\Rightarrow p\in]P_i,Q_i], \forall i\in[1,3] \\ m_i\in\mathcal{S}_i, \forall i\in[1,3] \\ \omega(m_1m_2m_3)\leq k-2}}\sum_{\substack{\frac{x}{m}<r\leq\frac{x+y_0}{m} \\ r\in\N_{\P}\cap\mathcal{S}_{\infty} \\ \omega(r)=k-\omega(m)}}1,\]
où $\P$ est l'ensemble des tous les nombres premiers sauf ceux appartenant à $\bigcup_{i=1}^3]P_i,Q_i]$. On définit $f:=\mathds{1}_{\N_{\P}}$ la fonction multiplicative indicatrice de l'ensemble $\N_{\P}$. On pose $k_m:=k-\omega(m)$, que l'on peut supposer supérieur ou égal à $2$, puisque $r\in\mathcal{S}_{\infty}$ et $mQ_{\infty}\leq Q_{\infty}\prod_{i=1}^3Q_i^{H_i\log Q_i}=X^{o(1)}$. On pose aussi $\kappa_m:=\frac{k_m-1}{\log_2 X}=\frac{k-1-\omega(m)}{\log_2 X}$. Ainsi, on a
\begin{equation*}
\sum_{\substack{\frac{x}{m}<r\leq\frac{x+y_0}{m} \\ r\in\N_{\P}\cap\mathcal{S}_{\infty} \\ \omega(r)=k-\omega(m)}}1=\Sigma+O\left(\Sigma_1+\Sigma_2+\Sigma_3\right),
\end{equation*}
où l'on a posé
\[\Sigma:=\sum_{\substack{\frac{x}{m}<r\leq\frac{x+y_0}{m} \\ r\in\mathcal{E}_{k_m}}}f(r),\]
et
\[\Sigma_1:=\sum_{\substack{\frac{x}{m}<r\leq\frac{x+y_0}{m} \\ r\in\mathcal{E}_{k_m} \\ \left(r,\prod_{P_{\infty}<p\leq Q_{\infty}}p\right)=1}}1,\hspace{.5cm}\Sigma_2:=\sum_{\substack{\frac{x}{m}<r\leq\frac{x+y_0}{m} \\ r\in\mathcal{E}_{k_m} \\ \exists p\in]P_{\infty},Q_{\infty}], p^2\vert r}}f(r),\hspace{.5cm}\Sigma_3:=\sum_{\substack{\frac{x}{m}<r\leq\frac{x+y_0}{m} \\ r\in\mathcal{E}_{k_m} \\ \exists p\in]P_{\infty},Q_{\infty}], \\ \exists \vert p-q\vert\leq\frac{2p}{H_{\infty}}, pq\vert r}}f(r).\]
Avec le Lemme~\ref{hyp1application}, on a facilement
\[\Sigma=\delta_{k_m}(X)\frac{y_0}{m}(1+o(1))\prod_{i=1}^{3}\left(\frac{\log P_i}{\log Q_i}\right)^{\kappa_m},\]
et $\Sigma_2+\Sigma_3\ll\left(\frac{1}{P_{\infty}}+\frac{\log\left(\frac{\log Q_{\infty}}{\log P_{\infty}}\right)}{H_{\infty}}\right)\Sigma=o(\Sigma)$. Par ailleurs, avec le Lemme~\ref{hyp2application}, il existe une constante $c=c(\varepsilon)>0$ fixée telle que $\Sigma_1\leq (1-c+o(1))\Sigma$. On obtient donc
\[\sum_{\substack{\frac{x}{m}<r\leq\frac{x+y_0}{m} \\ r\in\N_{\P}\cap\mathcal{S}_{\infty} \\ \omega(r)=k-\omega(m)}}1\asymp\delta_{k_m}(X)\frac{y_0}{m}\left(\frac{\log P_3}{\log Q_3}\right)^{\kappa_m}.\]
Par ailleurs, $\delta_{k_m}(X)\left(\frac{\log P_3}{\log Q_3}\right)^{\kappa_m}\asymp\delta_k(X)\e^{-\frac{k\log_3 X}{\log_2 X}}A^{\omega(m)}\prod_{a=1}^{\omega(m)}\left(1-\frac{a}{k}\right)$,
où l'on a posé $A:=\frac{k}{\log_2 X}\e^{\frac{\log_3 X}{\log_2 X}}$ Pour conclure, il suffit donc de montrer que l'on a
\[\mathcal{D}:=\sum_{\substack{m_i, \forall i\in[1,3] \\ p\vert m_i\Rightarrow p\in]P_i,Q_i], \forall i\in[1,3] \\ m_i\in\mathcal{S}_i, \forall i\in[1,3] \\ \omega(m_1m_2m_3)\leq k-2}}\frac{A^{\omega(m)}}{m}\prod_{a=1}^{\omega(m)}\left(1-\frac{a}{k}\right)\asymp\frac{k^2}{(\log_2 X)^2}\left(\e^{\frac{k\log_3 X}{\log_2 X}}-1\right).\]
On commence par minorer $\mathcal{D}$. En moyenne, un entier ayant $k$ facteurs premiers en a environ $\frac{k\log_3 X}{\log_2 X}\leq\log_3 X$ dans l'intervalle $]P_3,Q_3]$. On minore alors $\mathcal{D}$ efficacement en écrivant
\[\mathcal{D}\geq\sum_{\substack{p_1\in]P_1,Q_1] \\ p_2\in]P_2,Q_2]}}\frac{A^2}{p_1p_2}\sum_{N=1}^{\left\lceil\frac{2k\log_3 X}{\log_2 X}\right\rceil}\frac{A^N}{N!}\left(1-\frac{N+2}{k}\right)^{N+2}\sum_{q_1, ..., q_N\in]P_3, Q_3]}\frac{\mathds{1}_{q_1\cdots q_N\in\mathcal{S}_3}}{q_1\cdots q_N}.\]
On pose $M:=\sum_{P_3<p\leq Q_3}\frac{1}{p}=\log_3 X\left(1+O\left(\frac{1}{\log_2 X}\right)\right)$. On a alors, pour $2\leq N\leq\log_3 X$,
\[\sum_{q_1, ..., q_N\in]P_3, Q_3]}\frac{\mathds{1}_{q_1\cdots q_N\in\mathcal{S}_3}}{q_1\cdots q_N}=M^N+O\left(N^2\sum_{P_3<p\leq Q_3}\frac{M^{N-2}}{p^2}+N^2\sum_{\substack{P_3<p\leq Q_3 \\ \vert p-q\vert\leq\frac{2p}{H_3}}}\frac{M^{N-2}}{pq}\right)=M^N(1+o(1)).\]
Ainsi, puisque $AM=\frac{k\log_3 X}{\log_2 X}(1+o(1))$, on obtient
\[\mathcal{D}\gg A^2\sum_{N=1}^{\left\lceil\frac{2k\log_3 X}{\log_2 X}\right\rceil}\frac{(AM)^N}{N!}\gg A^2\left(\e^{AM}-1\right)\gg\frac{k^2}{(\log_2 X)^2}\left(\e^{\frac{k\log_3 X}{\log_2 X}}-1\right).\]
Il nous reste désormais à majorer $\mathcal{D}$. Pour cela, on écrit
\[\mathcal{D}\leq\prod_{i=1}^3\left\{\prod_{P_i<p_i\leq Q_i}\left(1+\frac{A}{p}\right)-1\right\},\]
ce qui fournit aisément le résultat recherché.

 \end{proof}

\begin{proof}[Démonstration du Théorème~$\ref{petitsk}$]
On peut déjà supposer $k\leq\frac{1}{4}\log_2 X$ grâce au Corollaire~\ref{corollaire}. Par ailleurs, le résultat est directement vrai lorsque $X<h\leq \delta_k(X)^{-1}X$ d'après le Lemme~$\ref{hyp1application}$. On pose désormais $\delta(X):=\frac{\delta_k(X)}{F_k(X)}$. En conservant les notations introduites depuis la sous-section~$\ref{sssec}$, d'après le lemme précédent, il suffit donc de démontrer que l'on a
\[\frac{1}{X}\int_X^{2X}\left\vert\frac{1}{h}\sum_{\substack{x<n\leq x+h \\ n\in\mathcal{E}_k\cap\mathcal{S}}}1-\frac{1}{y_0}\sum_{\substack{x<n\leq x+y_0 \\ n\in\mathcal{E}_k\cap\mathcal{S}}}1\right\vert^2\mathrm{d}x=o\left(\delta(X)^2\right)\]
dès que $\delta(X)^{-1}Q_1\leq h\leq X$ et $X$ tend vers l'infini. La preuve suit le même schéma que celle du Théorème~$\ref{th1}$. D'après les Lemmes~$\ref{Parseval}$ et~$\ref{MVT}$, il suffit de montrer que pour $T_0<T\leq X$, on a
\begin{equation}\label{equationintegrale}
\int_{T_0}^T\left\vert B(1+it)\right\vert^2\mathrm{d}t=o\left(\left(\frac{TQ_1}{\delta(X)X}+1\right)\delta(X)^2\right),
\end{equation}
lorsque $X$ tend vers l'infini, où l'on a posé
\[B(s):=\sum_{\substack{n\sim X\\ s\in\mathcal{E}_k\cap\mathcal{S}}}\frac{1}{n^s}\hspace{2cm}(s\in\C).\]
La factorisation $(\ref{facto})$ du polynôme de Dirichlet $B$ est toujours valable avec $\mathcal{A}=\mathcal{E}_k$ pour $\mathcal{A}'$ l'ensemble défini de la même manière que dans l'hypothèse~\ref{hyp5} du Théorème~$\ref{th1}$. On définit alors quatre ensembles $\mathcal{T}_1,\mathcal{T}_2,\mathcal{T}_3,\mathcal{U}$, de la même façon qu'à la section~$\ref{majorationdelintegrale}$ via~$(\ref{unionT0T})$ et~$(\ref{inegQ})$. On majore alors l'intégrale de~$(\ref{equationintegrale})$ successivement sur $\mathcal{T}_1$, puis $\mathcal{T}_2$ et $\mathcal{T}_3$, et enfin sur $\mathcal{U}$.\\

\textbf{Cas de $\mathcal{T}_1$ : } Avec~$(\ref{facto})$ et l'inégalité de Cauchy-Schwarz, on a
\[\int_{\mathcal{T}_1}\left\vert B(1+it)\right\vert^2\mathrm{d}t\ll\left\vert\mathcal{I}_1\right\vert\sum_{v\in\mathcal{I}_1}\int_{\mathcal{T}_1}\left\vert Q_{v,H_1}(1+it)R_{v,H_1}(1+it)\right\vert^2\mathrm{d}t+\int_{\mathcal{T}_1}\left\vert N_{H_1}(1+it)\right\vert^2\mathrm{d}t.\]
On traite d'abord l'intégrale de $\left\vert N_{H_1}(1+it)\right\vert^2$. Ce polynôme de Dirichlet est supporté sur les entiers $n\in\mathcal{E}_k\cap\mathcal{S}\cap\left]2X,2X\e^{1/H_1}\right]$. On majore son intégrale grâce au Lemme~$\ref{IMVT}$. Le Lemme~$\ref{sommepetitsk}$ permet de traiter la première somme associée à~$(\ref{eqIMVT})$. On majore maintenant la seconde somme. En décomposant tout entier $n\in\mathcal{E}_k\cap\mathcal{S}$ sous l'unique forme $m_1m_2m_3r$ où, pour $i\in\{1,3\}$, $m_i$ contient tous les facteurs premiers de $n$ dans $]P_i,Q_i]$, on obtient, pour $1\leq b\leq\frac{2X\e^{1/H_1}}{T}$,
\[\sum_{\substack{2X<n\leq 2X\e^{1/H_1} \\ n,n+b\in\mathcal{E}_k\cap\mathcal{S}}}1\leq\sum_{\substack{m_i,m_i',\ \forall 1\leq i\leq 3 \\ p\vert m_im_i'\Rightarrow p\in]P_i,Q_i],\ \forall 1\leq i\leq 3 \\ m_i,m_i'\in\mathcal{S}_i,\ \forall 1\leq i\leq 3 \\ \omega(m_1m_2m_3),\omega(m_1'm_2'm_3')\leq k-1 \\ (m_i,m_i')\vert b,\ \forall 1\leq i\leq 3}}\sum_{\substack{2X<n\leq 2X\e^{1/H_1} \\ m_1m_2m_3\vert n \\ m_1'm_2'm_3'\vert n+b \\ \omega_{\P}(n)= k-\omega(m_1m_2m_3) \\ \omega_{\P}(n+b)=k-\omega(m_1'm_2'm_3') \\ \frac{n}{m_1m_2m_3},\frac{n+b}{m_1'm_2'm_3'}\in\N_{\P}}}1,\]
où l'on pose encore $\P$ l'ensemble de tous les nombres premiers sauf ceux de $\bigcup_{i=1}^3]P_i,Q_i]$. Pour ce choix de $\P$, avec la notation~$(\ref{dEfEhetnP})$ et la formule de Mertens, on a $E(x)=\log_2 x\left(1-\frac{\log_3 x+O(1)}{\log_2 x}\right)$. On pose $m:=m_1m_2m_3$ et $m':=m_1'm_2'm_3'$. On remarque, dans la somme ci-dessus, que l'on a $m, m'\leq\prod_{i=1}^3Q_i^{H_i\log Q_i}=X^{o(1)}$. Pour majorer la dernière somme, on peut alors appliquer le point~\ref{deuxiemecascri} du Théorème~$\ref{nouvelleprop}$ en paramétrant $n$ modulo $\frac{mm'}{(m,m')}$. On obtient alors l'existence d'une constante $K>0$ telle que
\[\sum_{\substack{2X<n\leq 2X\e^{1/H_1} \\ m\vert n,\ m'\vert n+b \\ \omega_{\P}(n)=k-\omega(m) \\ \omega_{\P}(n+b)=k-\omega\left(m'\right) \\ \frac{n}{m},\frac{n+b}{m'}\in\N_{\P}}}1\ll\frac{b^K(m,m')}{\varphi(b)^K\varphi(m)\varphi(m')}\frac{X}{H_1(\log X)^2}\frac{E(X)^{2k-\omega(m)-\omega(m')-2}}{\left(k-\omega(m)-1\right)!\left(k-\omega(m')-1\right)!}.\]
Or, on a supposé $k\leq\frac{1}{4}\log_2 X\leq \frac{1}{2}E(X)$ lorsque $X$ est suffisamment grand. Donc, par un calcul élémentaire, on obtient
\[\frac{1}{(\log X)^2}\frac{E(X)^{2k-\omega(m)-\omega(m')-2}}{\left(k-\omega(m)-1\right)!\left(k-\omega(m')-1\right)!}\ll\delta_k(X)^2\e^{-\frac{2(k-1)\log_3 X}{\log_2 X}}\left(\frac{k}{E(X)}\right)^{\omega(m)+\omega(m')},\]
et ainsi
\begin{multline*}
\sum_{\substack{2X<n\leq 2X\e^{1/H_1} \\ n,n+b\in\mathcal{E}_k\cap\mathcal{S}}}1\ll\frac{b^K}{\varphi(b)^K}\delta_k(X)^2\frac{X}{H_1}\e^{-\frac{2k\log_3 X}{\log_2 X}}\\
\times\prod_{i=1}^3\sum_{\substack{m_i,m_i'\\ p\vert m_im_i'\Rightarrow p\in]P_i,Q_i]\\ m_i,m_i'\in\mathcal{S}_i\\ (m_i,m_i')\vert b}}\frac{(m_i,m_i')}{\varphi(m_i)\varphi(m_i')}\left(\frac{k}{E(X)}\right)^{\omega(m_i)+\omega(m_i')}.
\end{multline*}
Pour $i\in\{1,2\}$, on majore directement la somme ci-dessus en utilisant essentiellement le fait que $\log Q_i \ll\log P_i$.
\[\sum_{\substack{m_i,m_i'\\ p\vert m_im_i'\Rightarrow p\in]P_i,Q_i]\\ m_i,m_i'\in\mathcal{S}_i\\ (m_i,m_i')\vert b}}\frac{(m_i,m_i')}{\varphi(m_i)\varphi(m_i')}\left(\frac{k}{E(X)}\right)^{\omega(m_i)+\omega(m_i')}\leq\sum_{s,s'\geq 1}\left(\frac{k}{E(x)}\right)^{s+s'}\sum_{\substack{m_i,m_i' \\ p\vert m_im_i'\Rightarrow p\in]P_i,Q_i]\\ m_i,m_i'\in\mathcal{S}_i\\ (m_i,m_i')\vert b \\ \omega(m_i)=s,\ \omega(m_i')=s'}}\frac{(m_i,m_i')}{\varphi(m_i)\varphi(m_i')},\]
où la dernière somme se majore uniformément de manière grossière par
\[\prod_{P_i<p\leq Q_i}\left(1+\frac{2}{\varphi(p)}+\frac{p}{\varphi(p)^2}\right)\ll 1.\]
Il s'ensuit que pour $i\in\{1,2\}$, on a 
\[\sum_{\substack{m_i,m_i'\\ p\vert m_im_i'\Rightarrow p\in]P_i,Q_i]\\ m_i,m_i'\in\mathcal{S}_i\\ (m_i,m_i')\vert b}}\frac{(m_i,m_i')}{\varphi(m_i)\varphi(m_i')}\left(\frac{k}{E(X)}\right)^{\omega(m_i)+\omega(m_i')}\ll\left(\frac{k}{E(x)}\right)^2\ll\left(\frac{k}{\log_2 X}\right)^2.\]
Pour obtenir une majoration adéquate pour
\begin{equation*}
\sum_{1\leq b\leq\frac{2X\e^{1/H_1}}{T}}\sum_{\substack{2X<n\leq 2X\e^{1/H_1} \\ n,n+b\in\mathcal{E}_k\cap\mathcal{S}}}1,
\end{equation*}
il suffit alors de montrer que l'on a
\begin{equation*}
\e^{-\frac{2k\log_3 X}{\log_2 X}}\sum_{\substack{m_3,m_3'\\ p\vert m_3m_3'\Rightarrow p\in]P_3,Q_3]\\ m_3,m_3'\in\mathcal{S}_3\\ (m_3,m_3')\vert b}}\frac{(m_3,m_3')}{\varphi(m_3)\varphi(m_3')}\left(\frac{k}{E(X)}\right)^{\omega(m_3)+\omega(m_3')}\ll\frac{b}{\varphi(b)}\left(1-\e^{-\frac{k\log_3 X}{\log_2 X}}\right)^2.
\end{equation*}
On pose $A:=\frac{k}{E(X)}$, qui est inférieur à $\frac{1}{2}$ lorsque $X$ est suffisamment grand. La somme ci-dessus est alors inférieure à
\[\prod_{P_3<p\leq Q_3}\left(1+\frac{2A}{p-1}+\frac{pA^2\mathds{1}_{p\vert b}}{(p-1)^2}\right)-2\prod_{P_3<p\leq Q_3}\left(1+\frac{A}{p-1}\right)+1.\]
Cela est obtenu par multiplicativité, en ajoutant et retranchant les cas où l'un des $m_3, m_3'$ vaut~$1$, et en omettant le fait que $m_3$ et $m_3'$ n'ont pas deux facteurs premiers \og trop proches\fg. Avec des notations évidentes, on note la quantité ci-dessus $\Pi_1-2\Pi_2+1$. Par la formule de Mertens et avec la définition de $P_3, Q_3$, on a
\begin{align*}
\Pi_2&=\e^{A\log_3 X}\left(1+O\left(\frac{A}{\log_2 X}\right)\right),\\
\Pi_1&=\prod_{\substack{P_3<p\leq Q_3 \\ p\vert b}}\left(1+\frac{1}{p}\right)^{A^2}\e^{2A\log_3 X}\left(1+O\left(\frac{A}{\log_2 X}\right)\right).
\end{align*}
Lorsque $\frac{1}{2}\geq\frac{k}{\log_2 X}\geq \frac{1}{2\log_3 X}$, on a alors
\[\e^{-\frac{2k\log_3 X}{\log_2 X}}\left(\Pi_1-2\Pi_2+1\right)\ll\frac{b}{\varphi(b)}\ll\frac{b}{\varphi(b)}\left(1-\e^{-\frac{k\log_3 X}{\log_2 X}}\right)^2.\]
Et lorsque $\frac{1}{2\log_3 X}\geq \frac{k}{\log_2 X}\geq\frac{1}{\log_2 X}$, on a
\begin{align*}
\Pi_1-2\Pi_2+1&\ll\prod_{\substack{P_3<p\leq Q_3 \\ p\vert b}}\left(1+\frac{1}{p}\right)^{A^2}-1\\
&+\left(1+2A\log_3 X+O\left(A^2(\log_3 X)^2\right)\right)\left(1+O\left(\frac{A}{\log_2 X}\right)\right)\\
&-2\left(1+A\log_3 X+O\left(A^2(\log_3 X)^2\right)\right)\left(1+O\left(\frac{A}{\log_2 X}\right)\right)\\
&+1\\
&\ll A^2\log_3 X + A^2(\log_3 X)^2+\frac{A}{\log_2 X}\ll \frac{k^2(\log_3 X)^2}{(\log_2 X)^2}\ll\frac{b}{\varphi(b)}\left(1-\e^{-\frac{k\log_3 X}{\log_2 X}}\right)^2.
\end{align*}
Ainsi, on a la majoration
\begin{equation*}
\sum_{1\leq b\leq\frac{2X\e^{1/H_1}}{T}}\sum_{\substack{2X<n\leq 2X\e^{1/H_1} \\ n,n+b\in\mathcal{E}_k\cap\mathcal{S}}}1\ll\frac{X^2}{TH_1}\delta(X)^2,
\end{equation*}
qui est convenable. Pour majorer l'intégrale de $\left\vert Q_{v,H_1}(1+it)R_{v,H_1}(1+it)\right\vert^2$, il suffit d'appliquer la même méthode que dans la démonstration du Théorème~$\ref{th1}$ (\emph{cf.} page~\pageref{casj=1debut}), et d'utiliser le Théorème~\ref{nouvelleprop} de la même manière que ci-dessus. Il faut cependant décomposer les entiers $n\in\mathcal{A}'$ sous la forme $m_2m_3r$ seulement. Cela a pour incidence, par rapport au cas de $\int_{\mathcal{T}_1}\left\vert N_{H_1}(1+it)\right\vert^2\mathrm{d}t$ que l'on vient de traiter, de remplacer la quantité $\frac{\delta_k(X)}{F_k(X)}$ par $\frac{\delta_{k-1}(X)}{G_k(X)}$ où $G_k(X):=\frac{\log_2 X}{k}\left(1-\exp\left(-\frac{k\log_3 X}{\log_2 X}\right)\right)^{-1}$. Mais ces quantités sont en fait égales à un facteur borné près. On n'expose pas les calculs, qui sont plus simples que ceux du cas de $\mathcal{T}_2$, que l'on présente ci-dessous. On obtient alors
\[\int_{\mathcal{T}_1}\left\vert B(1+it)\right\vert^2\mathrm{d}t\ll\left(\frac{TQ_1}{\delta(X)X}+1\right)\delta(X)^2\left(H_1^2(\log Q_1)P_1^{-2\alpha_1}+\frac{1}{H_1}\right).\]

\textbf{Cas de $\mathcal{T}_2$ : }On applique la même méthode que dans la démonstration du Théorème~$\ref{th1}$, dont on reprend les notations (\emph{cf.} page~$\pageref{methodeth1}$). La majoration de $\int_{\mathcal{T}_2}\left\vert N_{H_2}(1+it)\right\vert^2\mathrm{d}t$ s'effectue exactement comme celle ci-dessus de $\int_{\mathcal{T}_1}\left\vert N_{H_1}(1+it)\right\vert^2\mathrm{d}t$. Pour le reste, on est amené à majorer, pour $u_1\in\mathcal{I}_1$ et $v_2\in\mathcal{I}_2$,
\[\int_{\mathcal{T}_{2, u_1}}\left\vert Q_{u_1,H_1}(1+it)^{\ell}R_{v_2,H_2}(1+it)\right\vert^2\mathrm{d}t,\]
où $\ell:=\left\lceil\frac{\log Y_2}{\log Y_1}\right\rceil=\left\lceil\frac{v_2/H_2}{u_{1}/H_{1}}\right\rceil$. On écrit le polynôme de Dirichlet $Q_{u_1,H_1}^{\ell}R_{v_2,H_2}$ sous la forme
\[Q_{u_1,H_1}(s)^{\ell}R_{v_2,H_2}(s)=\sum_{X\leq n\leq 2^{\ell+1}Y_1X}\frac{a_n}{n^s}\hspace{2cm}(s\in\C),\]
où les $a_n$ sont uniquement déterminés. Si $a_n$ est non nul, alors $n$ peut s'écrire sous la forme $n=m_1m_3\m r$ où $m_3\in\mathcal{S}_3$ est le produit des facteurs premiers de $n$ dans $]P_3,Q_3]$, $\m$ est un produit de $\ell$ facteurs premiers de $]Y_1,Y_1\e^{1/H_1}]\cap]P_1,Q_1]$, $m_1\in\mathcal{S}_1$ est le produit des facteurs premiers de $\frac{n}{\m}$ dans $]P_1,Q_1]$, $r\in\mathcal{S}_{\infty}$ n'a aucun facteur premier dans $]P_i,Q_i]$ pour $i\in\{1,3\}$, et $\omega(m_1m_3r)=k-1$. Ainsi, on a
\begin{equation}\label{majanll}
0\leq a_n\leq\sum_{\substack{m_1, m_3 \\ p\vert m_i\Rightarrow p\in]P_i,Q_i], \forall i\in\{1,3\} \\ m_i\in\mathcal{S}_i, \forall i\in\{1,3\} \\ \omega_{\P}(n)=k-1-\omega(m_1m_3)}}\sum_{\substack{p_1, ..., p_{\ell}\sim Y_1 \\ P_1<p_1, ..., p_{\ell}\leq Q_1 \\ m_1m_3p_1\cdots p_{\ell}\vert n \\ \frac{n}{m_1m_3p_1\cdots p_{\ell}}\in\N_{\P}\cap\mathcal{S}_{\infty}}}1\leq (\ell+1)!,
\end{equation}
où $\P$ est cette fois l'ensemble de tous les nombres premiers sauf ceux de $\bigcup_{i\in\{1,3\}}]P_i,Q_i]$. Comme précédemment, on majore alors l'intégrale de $\vert Q_{u_1,H_1}^{\ell}R_{v_2,H_2}\vert^2$ grâce au Lemme~\ref{IMVT}. On commence alors par traiter la somme associée à la première somme de~$(\ref{eqIMVT})$. On a
\begin{align*}
\sum_{X\leq n\leq2^{\ell+1}Y_1X}\left\vert\frac{a_n}{n}\right\vert^2&\leq\frac{(\ell+1)!}{X^2}\sum_{\substack{m_1, m_3 \\ p\vert m_i\Rightarrow p\in]P_i,Q_i], \forall i\in\{1,3\} \\ m_i\in\mathcal{S}_i, \forall i\in\{1,3\} \\ \omega(m_1m_3)\leq k-2}}\sum_{p_1, ..., p_{\ell}\sim Y_1}\sum_{\substack{\frac{X}{m_1m_3p_1\cdots p_{\ell}}\leq r\leq\frac{2^{\ell+1}Y_1X}{m_1m_3p_1\cdots p_{\ell}} \\ \omega(r)=k-1-\omega(m_1m_3) \\ r\in\N_{\P}}}1.
\end{align*}
On majore la dernière somme avec le point~\ref{premiercascri} du Théorème~\ref{nouvelleprop}. Ainsi, puisque $m_1m_3p_1\cdots p_{\ell}=X^{o(1)}$, avec $E(x)=\log_2 X\left(1-\frac{\log_3 X+O(1)}{\log_2 X}\right)$, on a
\begin{align*}
\sum_{\substack{\frac{X}{m_1m_3p_1\cdots p_{\ell}}\leq r\leq\frac{2^{\ell+1}Y_1X}{m_1m_3p_1\cdots p_{\ell}} \\ \omega(r)=k-1-\omega(m_1m_3) \\ r\in\N_{\P}}}1&\ll\frac{2^{\ell+1}Y_1X}{\log X}\frac{E(X)^{k-2-\omega(m_1m_3)}}{(k-2-\omega(m_1m_3))!}\frac{1}{m_1m_3p_1\cdots p_{\ell}}\\
&\ll \delta_{k-1}(X)2^{\ell}Y_1X\e^{-\frac{k\log_3 X}{\log_2 X}}\left(\frac{k}{E(X)}\right)^{\omega(m_1m_3)}\frac{1}{m_1m_3p_1\cdots p_{\ell}}
\end{align*}
En sommant sur les $p_1, ..., p_{\ell}$, puis $m_1$ et $m_3$, on obtient alors
\[\sum_{X\leq n\leq2^{\ell+1}Y_1X}\left\vert\frac{a_n}{n}\right\vert^2\ll\frac{(\ell+1)!2^{\ell}Y_1}{X}\delta(X).\]
Pour obtenir l'inégalité ci-dessus, la somme sur $m_3$ est la plus délicate à majorer, et se traite de manière analogue à la double somme sur $m_3, m_3'$ vue dans le cas de $\mathcal{T}_1$. On majore désormais la somme associée à la seconde somme de~$(\ref{eqIMVT})$. Pour cela, pour $1\leq b\leq\frac{2^{\ell+1}Y_1X}{T}$, en utilisant~$(\ref{majanll})$ et les décompositions $n=m_1m_3\mathfrak{m}r$ et $n+b=m_1'm_3'\mathfrak{m}'r'$, on écrit
\begin{multline}\label{niolkoplko}
\sum_{X\leq n\leq2^{\ell+1}Y_1X}\frac{\vert a_na_{n+b}\vert}{n(n+b)}\\
\ll\frac{1}{X^2}\sum_{\substack{m_1, m_1', m_3, m_3' \\ p\vert m_im_i'\Rightarrow p\in]P_i,Q_i], \forall i\in\{1,3\} \\ m_i, m_i'\in\mathcal{S}_i \\ \omega(m_1m_3), \omega(m_1'm_3')\leq k-2 \\ (m_i,m_i')\vert b,\forall i\in\{1,3\}}}\sum_{\substack{p_1, ..., p_{\ell}\sim Y_1 \\ q_1, ..., q_{\ell}\sim Y_1 \\ (m_1m_3p_1\cdots p_{\ell}, m_1'm_3'q_1\cdots q_{\ell})\vert b}}\sum_{\substack{X\leq n\leq 2^{\ell+1}Y_1X \\ m_1m_3p_1\cdots p_{\ell}\vert n \\ m_1'm_3'q_1\cdots q_{\ell}\vert n+b \\ \omega\left(\frac{n}{m_1m_3p_1\cdots p_{\ell}}\right)=k-1-\omega(m_1m_3) \\ \omega\left(\frac{n+b}{m_1'm_3'q_1\cdots q_{\ell}}\right)=k-1-\omega(m_1'm_3') \\ \frac{n}{m_1m_3p_1\cdots p_{\ell}}, \frac{n+b}{m_1'm_3'q_1\cdots q_{\ell}}\in\N_{\P}}}1.
\end{multline}
Grâce au deuxième point du Théorème~\ref{nouvelleprop}, toujours avec $E(X)=\log_2 X\left(1-\frac{\log_3 X+O(1)}{\log_2 X}\right)$, il existe une constante $K>0$ telle que la dernière somme ci-dessus soit
\begin{multline}\label{gvfgytg}
\ll \frac{b^K}{\varphi(b)^K}2^{\ell}Y_1X\delta_{k-1}(X)^2\e^{-\frac{2k\log_3 X}{\log_2 X}}\prod_{i\in\{1,3\}}\left(\frac{k}{E(X)}\right)^{\omega(m_i)+\omega(m_i')}\\
\times\frac{(m_1m_3p_1\cdots p_{\ell}, m_1'm_3'q_1\cdots q_{\ell})}{\varphi(m_1)\varphi(m_3)\varphi(p_1)\cdots \varphi(p_{\ell})\varphi(m_1')\varphi(m_3')\varphi(q_1)\cdots \varphi(q_{\ell})}.
\end{multline}
Par ailleurs, pour tout entier $g\geq 1$, on a
\[\sum_{\substack{1\leq b\leq\frac{2^{\ell+1}Y_1X}{T} \\ g\vert b}}\frac{b^K}{\varphi(b)^K}\ll\frac{g^{K-1}}{\varphi(g)^K}\frac{2^{\ell}Y_1X}{T}.\]
On majore alors
\[\sum_{1\leq b\leq\frac{2^{\ell+1}Y_1X}{T}}\sum_{X\leq n\leq2^{\ell+1}Y_1X}\frac{\vert a_na_{n+b}\vert}{n(n+b)}\]
en reportant~$(\ref{gvfgytg})$ dans~$(\ref{niolkoplko})$ et en sommant d'abord par rapport à la variable $b$. On est alors amené à majorer
\[\sum_{\substack{m_1, m_1', m_3, m_3' \\ p\vert m_im_i'\Rightarrow p\in]P_i,Q_i], \forall i\in\{1,3\} \\ m_i, m_i'\in\mathcal{S}_i}}\sum_{\substack{p_1, ..., p_{\ell}\sim Y_1 \\ q_1, ..., q_{\ell}\sim Y_1}}\frac{\frac{g^K}{\varphi(g)^K}\prod_{i\in\{1,3\}}\left(\frac{k}{E(X)}\right)^{\omega(m_i)+\omega(m_i')}}{\varphi(m_1)\varphi(m_3)\varphi(p_1)\cdots \varphi(p_{\ell})\varphi(m_1')\varphi(m_3')\varphi(q_1)\cdots \varphi(q_{\ell})}\]
de façon convenable. Par ailleurs, on a $\frac{g^K}{\varphi(g)^K}\leq 2^{2K\ell}\frac{(m_1, m_1')^K}{\varphi((m_1, m_1'))^K}\frac{(m_3, m_3')^K}{\varphi((m_3, m_3'))^K}$. Ainsi, en suivant les calculs effectués dans le cas de $\mathcal{T}_1$, sachant que $\sum_{p\sim Y_1}\frac{1}{\varphi(p)}\leq 2$, on trouve facilement
\begin{equation*}\sum_{1\leq b\leq\frac{2^{\ell+1}Y_1X}{T}}\sum_{X\leq n\leq2^{\ell+1}Y_1X}\frac{\vert a_na_{n+b}\vert}{n(n+b)}\ll\frac{2^{(2K+4)\ell}Y_1^2\delta(X)^2}{T},
\end{equation*}
et donc
\begin{equation*}
\int_{\mathcal{T}_{2, u_1}}\left\vert Q_{u_1,H_1}(1+it)^{\ell}R_{v_2,H_2}(1+it)\right\vert^2\mathrm{d}t\ll\left(\frac{T}{\delta(X)X}+1\right)\delta(X)^2Y_1^2(\ell !)^{1+o(1)}.
\end{equation*}
Finalement, on obtient
\[\int_{\mathcal{T}_2}\vert B(1+it)\vert^2\mathrm{d}t\ll\left(\frac{T}{\delta(X)X}+1\right)\delta(X)^2\left(\vert\mathcal{I}_2\vert^2\vert\mathcal{I}_1\vert Q_1^{2+2\alpha_1}(\ell !)^{1+o(1)}Y_2^{-2(\alpha_2-\alpha_1)}+\frac{1}{H_2}\right).\]
On a par ailleurs $\ell\log\ell\leq(\log Y_2)\frac{\log_2 Q_2}{\log P_1-1}+\log_2 Q_2 +1$. Ainsi, avec~$(\ref{alphapetitsk})$, le terme entre parenthèses après $\delta(X)^2$ ci-dessus tend bien vers $0$ lorsque $X$ tend vers l'infini.

\textbf{Cas de $\mathcal{T}_3$ : }L'intégrale sur $\mathcal{T}_3$ se traite de la même manière que celle sur $\mathcal{T}_2$, \emph{mutatis mutandis}.

\textbf{Cas de $\mathcal{U}$ : }Il suffit de suivre le cas de $\mathcal{U}$ dans la preuve du Théorème~$\ref{th1}$, \emph{mutatis mutandis}.
\end{proof}

\bibliographystyle{plain-fr}
\bibliography{Bibliographie}

\noindent\small\bsc{École Normale Supérieure, 45 rue d’Ulm 75230 Paris Cedex 05, France}\\

\noindent\bsc{Institut de Math\'ematiques de Jussieu-PRG, Universit\'e Paris Diderot,
Sorbonne Paris Cit\'e, 75013 Paris, France}\\

\noindent\textit{E-mail :} \url{eliegoudout@hotmail.com}

\end{document}